\documentclass[onecollarge,referee,12pt]{svjour2}       
\smartqed  
\usepackage[letterpaper,text={6.5in,8.6in},centering]{geometry}
\usepackage{amssymb,amsmath,times,graphicx,multirow,subfigure,url,warmread}
\usepackage[all,import]{xy}

\newcommand{\norm}[1]{\ensuremath{\left\| #1 \right\|}}
\newcommand{\bracket}[1]{\ensuremath{\left[ #1 \right]}}
\newcommand{\braces}[1]{\ensuremath{\left\{ #1 \right\}}}
\newcommand{\parenth}[1]{\ensuremath{\left( #1 \right)}}
\newcommand{\refeqn}[1]{(\ref{eqn:#1})}
\newcommand{\reffig}[1]{Fig. \ref{fig:#1}}
\newcommand{\tr}[1]{\mbox{tr}\ensuremath{\negthickspace\bracket{#1}}}
\newcommand{\deriv}[2]{\ensuremath{\frac{\partial #1}{\partial #2}}}
\newcommand{\SO}{\ensuremath{\mathrm{SO(3)}}}
\newcommand{\T}{\ensuremath{\mathrm{T}}}
\newcommand{\so}{\ensuremath{\mathfrak{so}(3)}}
\newcommand{\SE}{\ensuremath{\mathrm{SE(3)}}}
\renewcommand{\Re}{\ensuremath{\mathbb{R}}}
\def\makeheadbox{{%
\hbox to0pt{\vbox{\baselineskip=10dd\hrule\hbox
to\hsize{\vrule\kern3pt\vbox{\kern3pt \hbox{\bfseries Journal of
Optimization Theory and Applications manuscript No.} \hbox{  }
\kern3pt}\hfil\kern3pt\vrule}\hrule}%
\hss}}}

\journalname{Journal of Optimization Theory and Applications}

\begin{document}

\title{Optimal Attitude Control of a Rigid Body\\ using Geometrically Exact Computations on \SO
\thanks{The research of ML has been supported in part by NSF under grant DMS-0504747, and by a grant from the Rackham Graduate School, University of Michigan.
The research of NHM has been supported in part by NSF under grant
ECS-0244977.} }

\author{Taeyoung Lee\and
        Melvin Leok\and
        N. Harris McClamroch
}

\institute{Taeyoung Lee, Graduate Student \at
           Department of Aerospace Engineering, University of
           Michigan, Ann Arbor, MI 48109.
           \email{tylee@umich.edu}
        \and
           Melvin Leok, T. H. Hildebrandt Research Assistant Professor,  \at
           Department of Mathematics, University of Michigan, Ann
           Arbor, MI 48109.
           \email{mleok@umich.edu}
        \and
           N. Harris McClamroch, Professor,\at
           Department of Aerospace Engineering, University of
           Michigan, Ann Arbor, MI 48109.
           \email{nhm@engin.umich.edu}
}

\date{January 11, 2006}

\maketitle

\begin{abstract}
An efficient and accurate computational approach is proposed for
optimal attitude control of a rigid body.  The problem is formulated
directly as a discrete time optimization problem using a Lie group
variational integrator.  Discrete necessary conditions for
optimality are derived, and an efficient computational approach is
proposed to solve the resulting two point boundary value problem.
The use of geometrically exact computations on $\SO$ guarantees that
this optimal control approach has excellent convergence properties
even for highly nonlinear large angle attitude maneuvers. Numerical
results are presented for attitude maneuvers of a 3D pendulum and a
spacecraft in a circular orbit. 
\end{abstract}

\section{Introduction}
A discrete optimal control problem for attitude dynamics of a rigid
body in the presence of an attitude dependent potential is
considered. The objective is to minimize the square of the $l_2$
norm of external control torques which transfer a given initial
attitude and an initial angular momentum of the rigid body to a
desired terminal attitude and a terminal angular momentum during a
fixed maneuver time. The attitude of the rigid body is defined by
the orientation of a body fixed frame with respect to a reference
frame; the attitude is represented by a rotation matrix that is a
$3\times 3$ orthogonal matrix with determinant of 1. Rotation
matrices have a Lie group structure denoted by $\SO$.

The dynamics of a rigid body has fundamental invariance properties.
In the absence of nonconservative forces, total energy is preserved.
A consequence of Noether's theorem is that symmetries in the
Lagrangian result in conservation of the associated momentum map.
Furthermore, the configuration space of the rigid body has the
orthogonal structure of the Lie group $\SO$. General-purpose
numerical integration methods, including the popular Runge--Kutta
schemes, typically preserve neither first integrals nor the
geometric  characteristics of the configuration space. In
particular, the orthogonal structure of the rotation matrices is not
preserved numerically with standard schemes.

A Lie group variational integrator that preserves those geometric
features is presented in~\cite{Mel.Pre04}, and integrators on the
configuration space $\SO$ and $\SE$ are developed
in~\cite{pro:cca05} and  \cite{jo:CMAME05}, respectively. These
integrators are obtained from a discrete variational principle, and
they exhibit the characteristic symplectic and momentum preservation
properties, and good energy behavior characteristic of variational
integrators~\cite{jo:marsden}. Since the rotation matrices are
updated by a group operation, they automatically evolve on the
rotation group without the need for reprojection techniques or
constraints~\cite{ic:iserles}.

Optimal attitude control problems are studied in~\cite{Spin.MCSS98}.
The angular velocity of a rigid body is treated as a control input;
an optimal angular velocity that steers the rigid body is derived
from the attitude kinematics. Continuous time optimal control
problems on a Riemannian manifold are studied in~\cite{IsBl.CDC04},
where necessary conditions for optimality are derived from a
variational principle. An optimal control problem based on discrete
mechanics is studied in~\cite{JunMarObe.IFAC05}. The discrete
equations of motion and the boundary conditions are imposed as
constraints, and the optimal control problem is solved by a
general-purpose parameter optimization tool. This approach requires
large computation time, since the number of optimization parameters
is proportional to the number of discrete time steps. Discrete
necessary conditions for optimal control of the attitude dynamics of
a rigid body are presented in~\cite{HusLeoSanBlo.ACC06}.

This paper proposes an exact and efficient computational approach to
solve an optimal control problem associated with the attitude
dynamics of a rigid body that evolves on the configuration space
$\SO$. We assume that the control inputs are parameterized by their
value at each time step. A Lie group variational integrator on $\SO$
that includes the effects of external control inputs is developed
using the discrete Lagrange-d'Alembert principle. Discrete necessary
conditions for optimality are obtained using a variational
principle, while imposing the Lie group variational integrator as
dynamic constraints.

The necessary conditions are expressed as a two point boundary value
problem on $\T^*\SO$ and its dual. Sensitivity derivatives along an
extremal solution are developed by following the procedures
presented in~\cite{pro:acc06}, and they are used to construct an
algorithm that solves the boundary value problem efficiently. Since
the attitude of the rigid body is represented by a rotation matrix,
and the orthogonal structure of rotation matrices is preserved by
the Lie group variational integrator, the discretization of the
optimal control problem does not exhibit singularities.

This paper is organized as follows. In Section \ref{sec:eom}, a Lie
group variational integrator is developed using the discrete
Lagrange-d'Alembert principle. Necessary conditions for optimality
and a proposed approach to solve the two point boundary problem are
presented in Section \ref{sec:opt}. Numerical results for the
attitude control of an underactuated 3D pendulum, and for a fully
actuated spacecraft in a circular orbit are given in Section
\ref{sec:sim}.

\section{Equations of Motion for the Attitude Dynamics of a Rigid
Body}\label{sec:eom}

In this section, we define a rigid body model in a potential field
and we develop discrete equations of motion for the attitude
dynamics of the rigid body, referred to as a Lie group variational
integrator. These discrete equations of motion are used as dynamic
constraints for the optimal control problem presented in Section
\ref{sec:opt}.

\subsection{Rigid body model}
We consider the attitude dynamics of a rigid body in the presence of
an attitude dependent potential. The configuration space is the Lie
group, $\SO$. We assume that the potential $U(\cdot):\SO\mapsto\Re$
is determined by the attitude of the rigid body, $R\in\SO$. External
control inputs generate moments about the mass center of the rigid
body. A spacecraft on a circular orbit including gravity gradient
effects~\cite{pro:acc06}, a 3D pendulum~\cite{pro:cca05}, or a free
rigid body can be modeled in this way. The continuous equations of
motion are given by
\begin{gather}
\dot\Pi + \Omega\times \Pi = M+Bu,\\
\dot{R} = R S(\Omega),\label{eqn:Rdot}
\end{gather}
where $\Omega\in\Re^3$ is the angular velocity of the body expressed
in the body fixed frame, and $\Pi=J\Omega\in\Re^3$ is the angular
momentum of the body for a moment of inertia matrix
$J\in\Re^{3\times 3}$. $M\in\Re^3$ is the moment due to the
potential, $u\in\Re^m$ is the external control input, and
$B\in\Re^{3\times m}$ is an input matrix. If the rank of the input
matrix is less than 3, then the rigid body is underactuated. The
matrix valued function, $S(\cdot):\Re^3\mapsto \so$ is a skew
mapping defined such that $S(x)y=x\times y$ for all $x,y\in\Re^3$.
The Lie algebra $\so$ is identified by $3\times 3$  The moment due
to the potential is determined by the relationships,
$S(M)=\deriv{U}{R}^TR-R^T\deriv{U}{R}$, or more explicitly,
\begin{gather}
M=r_1\times v_{r_1} + r_2\times v_{r_2} +r_3\times
v_{r_3},\label{eqn:M}
\end{gather}
where $r_i,v_{r_i}\in\Re^{1\times 3}$ are the $i$th row vectors of
$R$ and $\deriv{U}{R}$, respectively. A detailed derivation of the
above equations can be found in~\cite{pro:cca05}.

\subsection{Lie group variational integrator}
The attitude dynamics of a rigid body exhibit geometric invariant
features. In the absence of an external control input, the total
energy is preserved. If there is a symmetry in the potential
function, the corresponding momentum map is preserved. The attitude
as described by a rotation matrix is always orthogonal. Classical
numerical integration methods typically preserve neither first
integrals nor the geometry of the configuration space, $\SO$. In
particular, standard Runge-Kutta method fail to capture the energy
dissipation of a controlled system accurately~\cite{jo:marsden}.

It is often proposed to parameterize \refeqn{Rdot} by Euler angles
or quaternions instead of integrating \refeqn{Rdot} directly.
However, Euler angles have singularities, and the unit length of a
quaternion vector is not preserved by classical numerical
integration. Furthermore, renormalizing the quaternion vector at
each step tends to break other conservation properties.

We describe a Lie group variational integrator that respects these
geometric properties. It is obtained from a discrete variational
approach, and therefore it exactly preserves the momentum and
symplectic form, while exhibiting good energy behavior over
exponentially long times. Since a Lie group numerical
method~\cite{ic:iserles} is explicitly adopted, the rotation matrix
automatically remains on $\SO$.

The Lie group variational integrator is obtained by following
procedures commonly adopted in Lagrangian mechanics. The variational
approach is based on discretizing Hamilton's principle rather than
discretizing the continuous equations of motion. The velocity phase
space of the continuous Lagrangian is replaced by discrete
variables, and a discrete Lagrangian is chosen. Taking a variation
of the action sum defined as the summation of the discrete
Lagrangian, we obtain a Lagrangian form of the discrete equations of
motion using the Lagrange-d'Alembert principle. A discrete version
of the Legendre transformation yields a Hamiltonian form.

The detailed derivation is presented in~\cite{pro:cca05}
and~\cite{jo:CMAME05}. In this paper, we extend these results to
include the effects of external control inputs. Consider the fixed
integration step size $h\in\Re$. Let $R_k \in \SO$ denote the
attitude of the rigid body at time $t=kh$. We introduce a new
variable $F_k\in\SO$ defined by
\begin{align}
F_k = R_k^T R_{k+1},
\end{align}
which represents a relative attitude between integration steps. If
we find $F_k\in\SO$ at each integration step, the rotation matrix is
updated by multiplication of two rotation matrices, i.e.
$R_{k+1}=R_kF_k$, which is a group operation on $\SO$. This
guarantees that the rotation matrix evolves on $\SO$ automatically.
This is the approach of Lie group methods~\cite{ic:iserles}. The
following procedure provides an expression for $F_k$ using the
discrete Lagrange-d'Alembert principle.

Using the kinematic relationship \refeqn{Rdot}, $S(\Omega_k)$ can be
approximated as
\begin{align*}
S(\Omega_k) = R_k^T\dot{R}_k \approx R_k^T
\parenth{\frac{R_{k+1}-R_k}{h}} =
\frac{1}{h}\parenth{F_k-I_{3\times3}}.
\end{align*}
Using the above equation, we can show that the kinetic energy of the
rigid body is given by
\begin{align*}
T&=\frac{1}{2} \tr{S(\Omega_k)J_dS(\Omega_k)^T}= \frac{1}{h^2}
\tr{\parenth{I_{3\times3}-F_k}J_d},
\end{align*}
where $J_d\in\Re^{3\times 3}$ is a nonstandard moment of inertia
matrix of the rigid body defined in terms of the standard moment of
inertia matrix $J\in\Re^{3\times 3}$ as
$J_d=\frac{1}{2}\tr{J}I_{3\times 3}-J$. Define a discrete Lagrangian
$L_d$ as
\begin{align}
L_d (R_k, F_k) & = \frac{1}{h} \tr{\parenth{I_{3\times3}-F_k}J_d} -
h U(R_{k+1}).\label{eqn:Ld0}
\end{align}
This discrete Lagrangian is a first-order approximation of the
integral of the continuous Lagrangian over one integration step.
Therefore, the following action sum, defined as the summation of the
discrete Lagrangian, approximates the action integral;
$\mathfrak{G}_d = \sum_{k=0}^{N-1}L_d(R_k,F_k)$.
Taking a variation of the action sum, we obtain the discrete
equations of motion using the discrete Lagrange-d'Alembert
principle. The variation of a rotation matrix can be expressed using
the exponential of a Lie algebra element:
\begin{align*}
R_k^\epsilon = R_k e^{\epsilon\eta_k},
\end{align*}
where $\epsilon \in \Re$ and $\eta_k \in \so$ is the variation
expressed as a skew symmetric matrix. Thus the infinitesimal
variation is given by
\begin{align}
\delta R_k & = \frac{d}{d\epsilon} \bigg|_{\epsilon=0} R_k^\epsilon
= R_k \eta_k.\label{eqn:delRk0}
\end{align}

The Lagrange-d'Alembert principle states that the following equation
is satisfied for all possible variations $\eta_k\in\so$.
\begin{align}
\delta\sum_{k=0}^{N-1} \frac{1}{h}
\tr{\parenth{I_{3\times3}-F_k}J_d} - h
U(R_{k+1})-\sum_{k=0}^{N-1}\frac{h}{2}\tr{\eta_{k+1}S(Bu_{k+1})}=0.
\end{align}
Using the expression of the infinitesimal variation of a rotation
matrix \refeqn{delRk0} and using the fact that the variations vanish
at the end points, the above equation can be written as
\begin{align*}
\sum_{k=1}^{N-1} \tr{\eta_k \braces{\frac{1}{h} \parenth{F_k J_d
-J_dF_{k-1}} + h R_k^T\deriv{U}{R_k}- \frac{h}{2}S(Bu_k)}}=0.
\end{align*}
Since the above expression should be zero for all possible
variations $\eta_k\in\so$, the expression in the braces should be
symmetric. Then, \textit{the discrete equations of motion in
Lagrangian form} are given by
\begin{gather}
\frac{1}{h} \parenth{F_{k+1} J_d -J_dF_{k}-J_d F_{k+1}^T+ F_{k}^T
J_d }  = h
S(M_{k+1})+hS(Bu_{k+1}),\label{eqn:updatef}\\
R_{k+1}=R_k F_k\label{eqn:updateR0}.
\end{gather}
Using the discrete version of the Legendre transformation,
\textit{the discrete equations of motion in Hamiltonian form} are
given by
\begin{gather}
h S(\Pi_k) = F_k J_d - J_dF_k^T,\label{eqn:findf}\\
R_{k+1} = R_k F_k,\label{eqn:updateR}\\
\Pi_{k+1} = F_k^T \Pi_k + h
\parenth{M_{k+1}+Bu_{k+1}}.\label{eqn:updatePi}
\end{gather}
Given $(R_k,\Pi_k)$, we can obtain $F_k$ by solving \refeqn{findf},
and $R_{k+1}$ is obtained by \refeqn{updateR}. The moment due to the
potential $M_{k+1}$ can be calculated  by \refeqn{M}. Finally,
$\Pi_{k+1}$ is updated by \refeqn{updatePi}. This yields a map
$(R_k, \Pi_k)\mapsto(R_{k+1},\Pi_{k+1})$, and this process can be
repeated. The only implicit part is solving \refeqn{findf}. We can
express \refeqn{findf} in terms of a Lie algebra element
$S(f_k)=\mathrm{logm}(F_k)\in\so$, and find $f_k\in\Re^3$
numerically by a Newton iteration. The relative attitude $F_k$ is
obtained by the exponential map: $F_k=e^{S(f_k)}$. Therefore we are
guaranteed that $F_k$ is a rotation matrix.

The order of the variational integrator is equal to the order of the
corresponding discrete Lagrangian. Consequently, the above Lie group
variational integrator is of first order since \refeqn{Ld0} is a
first-order approximation. While higher-order variational
integrators can be obtained by modifying \refeqn{Ld0}, we use the
first-order integrator because it yields a compact form for the
necessary conditions that preserves the geometry; these necessary
conditions are developed in Section \ref{sec:opt}.

\section{Discrete Optimal Control of the Attitude Dynamics of a Rigid
Body}\label{sec:opt} We formulate a discrete optimal control problem
for the attitude dynamics of a rigid body, and we derive necessary
conditions for optimality using a variational principle. The
necessary conditions are expressed as a two point boundary value
problem, and a computational approach to solve the boundary value
problem is proposed using sensitivity derivatives.

\subsection{Problem formulation}
A discrete time optimal control problem in $\SO$ is formulated as a
maneuver of a rigid body from a given initial attitude $R_0\in\SO$
and an initial angular momentum $\Pi_0\in\Re^3$ to a desired
terminal attitude $R_N^d\in\SO$ and a terminal angular momentum
$\Pi_N^d\in\Re^3$ during a given maneuver time $N$. The performance
index is the square of the $l_2$ norm of the control input; the
discrete equations of motion developed in the previous section are
imposed as constraints.
\begin{gather*}
\text{given: } R_0,\Pi_0,\,R_N^d,\Pi_N^d,\,N,\\
\min_{u_{k+1}} \mathcal{J}=\sum_{k=0}^{N-1} \frac{h}{2}\norm{u_{k+1}}^2,\\
\text{such that } R_N=R_N^d,\,\Pi_N=\Pi_N^d,\\ \text{  subject to
\refeqn{findf}, \refeqn{updateR} and \refeqn{updatePi}.}
\end{gather*}

In~\cite{JunMarObe.IFAC05}, an optimal control problem based on
discrete mechanics is considered. The control inputs at each
discrete step are considered as optimization parameters, and the
discrete equations of motion and the boundary conditions are imposed
as constraints. The optimization problem is solved numerically by a
general-purpose parameter optimization tool such as Sequential
Quadratic Programming (SQP). The same approach can be applied to the
above optimization problem. However, it has a large computational
burden since the number of optimization parameters, $m\times N$, is
proportional to the number of integration steps. Usually, a large
time step size is chosen to make the number of integration steps
small, or the control inputs are approximated by collocation points.
The resulting control inputs tends to be under-resolved and
sub-optimal.

We derive necessary conditions for optimality using the standard
calculus of variations. We assume that the control inputs are
parameterized by their value at each time step. The necessary
conditions are expressed as a two point boundary value problem.

\subsection{Necessary conditions of optimality}

Define an augmented performance index as
\begin{align}
\mathcal{J}_a = \sum_{k=0}^{N-1} & \frac{h}{2}\norm{u_{k+1}}^2
+\lambda_k^{1,T}S^{-1}\!\parenth{\mathrm{logm}(F_k-R_{k}^TR_{k+1})}\nonumber\\
&\quad +\lambda_k^{2,T}\braces{-\Pi_{k+1} + F_k^T \Pi_k +
h\parenth{M_{k+1}+Bu_{k+1}}},\label{eqn:Ja}
\end{align}
where $\lambda_k^{1}, \lambda_k^{2}\in \Re^3$, are Lagrange
multipliers corresponding to the discrete equations of motion. The
augmented performance index is chosen such that the dimension of the
multipliers is equal to the dimension of the rotation matrix and the
angular momentum vector. The discrete kinematics equation
\refeqn{updateR} is transformed into a matrix logarithm form. The
constraints arising from the discrete kinematics equation
\refeqn{updateR} and the angular momentum equation \refeqn{updatePi}
are explicitly applied. Equation \refeqn{findf} appears in the
discrete equations of motion because we introduce the auxiliary
variable $F_k\in\SO$. The constraint \refeqn{findf} is considered
implicitly when taking a variation of the performance index.

Consider small variations from a given trajectory denoted by
$\Pi_k,R_k,F_k,u_k$:
\begin{align}
\Pi_k^\epsilon & = \Pi_k + \epsilon \delta \Pi_k,\label{eqn:delPik}\\
R_k^\epsilon & = R_k e^{\epsilon S(\zeta_k)},\nonumber\\
& = R_k + \epsilon R_k S(\zeta_k) + \mathcal{O}(\epsilon^2),\label{eqn:delRk}\\
F_k^\epsilon & = F_k e^{\epsilon S(\xi_k)},\nonumber\\
& = F_k + \epsilon F_k S(\xi_k) +
\mathcal{O}(\epsilon^2),\label{eqn:delFk}
\end{align}
where $\zeta_k,\xi_k\in\Re^3\simeq\so$. The real space $\Re^3$ is
isomorphic to the Lie algebra $\so$ according to the skew mapping
$S(\cdot):\Re^3\mapsto\so$. The variations of the rotation matrices
are expressed using the exponential of the Lie algebra elements. The
corresponding infinitesimal variations of $\Pi_k,R_k$, and $F_k$ are
given by $\delta\Pi_k,\,\delta R_k=R_kS(\zeta_k)$, and $\delta F_k=
F_kS(\xi_k)$, respectively.

The variation of the augmented performance index is obtained from
the above expressions. Instead of taking a variation of the matrix
logarithm in \refeqn{Ja}, we take a variation of the kinematics
equation, \refeqn{updateR} and we use it as a constrained variation.
Since $F_k=R_k^T R_{k+1}$ by \refeqn{updateR}, the variation $\delta
F_k$ is given by
\begin{gather*}
\delta F_k = \delta R_k^T R_{k+1} + R_k^T \delta R_{k+1}.
\end{gather*}
Substituting the expression for the infinitesimal variation of
$R_k$, we obtain
\begin{gather*}
 F_kS(\xi_k) = -S(\zeta_k)F_k + F_k S(\zeta_{k+1}).
\end{gather*}
Multiplying both sides of the above equation by $F_k^T$ and using
the property $S(R^T x)=R^TS(x)R$ for all $R\in\SO$ and $x\in\Re^3$,
we obtain
\begin{align}
\xi_k = -F_k^T \zeta_k +\zeta_{k+1}.\label{eqn:xik0}
\end{align}
We use \refeqn{xik0} as a constrained variation equivalent to
\refeqn{updateR}.

Now we develop another expression for a constrained variation using
\refeqn{findf}. Since we do not use \refeqn{findf} explicitly as a
constraint in \refeqn{Ja}, $\delta\Pi_k$ and $\delta F_k$ are not
independent. Taking a variation of \refeqn{findf}, we obtain
\begin{align*}
h S(\delta\Pi_k) & = F_k S(\xi_k) J_d +  J_d S(\xi_k)F_k^T.
\end{align*}
Using the properties, $S(Rx)=RS(x)R^T$ and
$S(x)A+A^TS(x)=S(\braces{\tr{A}I_{3\times 3}-A}x)$ for all
$x\in\Re^3$, $A\in\Re^{3\times 3}$, and $R\in\SO$, the above
equation can be written as
\begin{align*}
h S(\delta\Pi_k) & = S(F_k\xi_k) F_kJ_d +  J_d F_k^T S(F_k \xi_k),\\
& = S(\braces{\tr{F_kJ_d}I_{3\times 3}-F_kJ_d}F_k\xi_k).
\end{align*}
Thus, $\xi_k$ is given by
\begin{align}
\xi_k = \mathcal{B}_k \delta\Pi_k, \label{eqn:xik}
\end{align}
where $\mathcal{B}_k=hF_k^T\braces{\tr{F_kJ_d}I_{3\times
3}-F_kJ_d}^{-1}\in\Re^{3\times 3}$. Equation \refeqn{xik} shows the
relationship between $\delta\Pi_k$ and $\delta F_k$.

Since the moment due to the potential $M_k$ is dependent on the
attitude of the rigid body, the variation of the moment $\delta M_k$
can be written using a variation of the rotation matrix:
\begin{align}
\delta M_k =\mathcal{M}_k\zeta_k,\label{eqn:delMk}
\end{align}
where $\mathcal{M}_k\in\Re^{3\times 3}$ is expressed in terms of the
the attitude of the rigid body, and the expression is determined by
the potential field. We present detailed expressions of
$\mathcal{M}_k$ in Section \ref{sec:sim} for a 3D pendulum and for a
spacecraft in a circular orbit. Using \refeqn{xik0} and
\refeqn{xik}, $\delta M_{k+1}$ is given by
\begin{align}
\delta M_{k+1} & = \mathcal{M}_{k+1} \zeta_{k+1},\nonumber\\
& = \mathcal{M}_{k+1}F_k^T \zeta_k
+\mathcal{M}_{k+1}\mathcal{B}_k\delta\Pi_k.\label{eqn:delMkp}
\end{align}

Now, we take a variation of the augmented performance index
\refeqn{Ja} using the constrained variations \refeqn{xik0},
\refeqn{xik}, and \refeqn{delMkp}. Using \refeqn{xik0}, the
variation of the performance index is given by
\begin{align}
\delta\mathcal{J}_a &= \sum_{k=0}^{N-1} h\delta u_{k+1}^T u_{k+1}
+ \lambda_k^{1,T}\braces{\xi_k +F_k^T \zeta_k -\zeta_{k+1}}\nonumber\\
& \quad\quad\quad+\lambda_k^{2,T}\braces{-\delta \Pi_{k+1} + \delta
F_k^T \Pi_k+F_k^T \delta\Pi_k + h \delta M_{k+1}+hB\delta u_{k+1}}
.\label{eqn:delJa0}
\end{align}
Substituting \refeqn{delMkp} into \refeqn{delJa0} and rearranging,
we obtain
\begin{align}
\delta\mathcal{J}_a &= \sum_{k=0}^{N-1} h\delta u_{k+1}^T
\braces{u_{k+1}+B^T\lambda_k^2}
-\zeta_{k+1}^T\lambda_k^1+\zeta_k^T \braces{F_k \lambda_k^1+hF_k \mathcal{M}_{k+1}^T\lambda_k^2}\nonumber\\
&\quad\quad\quad-\delta\Pi_{k+1}^T \lambda_k^2+\delta\Pi_k^T \braces{F_k\lambda_k^2+h\mathcal{B}_k^T \mathcal{M}_{k+1}^T\lambda_k^2}\nonumber\\
&\quad\quad\quad +\xi_k^T
\braces{-S(F_k^T\Pi_k)\lambda_k^2+\lambda_k^1}.\label{eqn:delJa01}
\end{align}
Substituting \refeqn{xik} into \refeqn{delJa01} and using the fact
that the variations $\zeta_k,\,\delta\Pi_k$ vanish at $k={0,N}$, we
obtain
\begin{align}
\delta\mathcal{J}_a &= \sum_{k=1}^{N-1} h\delta u_{k}^T
\braces{u_{k}+B^T\lambda_{k-1}^2}
+\zeta_k^T \braces{-\lambda_{k-1}^1+F_k \lambda_k^1+hF_k \mathcal{M}_{k+1}^T\lambda_k^2}\nonumber\\
&\quad\quad\quad+\delta\Pi_k^T
\braces{-\lambda_{k-1}^2+(F_k-\mathcal{B}_k^T
S(F_k^T\Pi_k)+h\mathcal{B}_k^T
\mathcal{M}_{k+1}^T)\lambda_k^2+\mathcal{B}_k^T\lambda_k^1}.
\end{align}

Since $\delta\mathcal{J}_a=0$ for all variations of $\delta u_k,
\zeta_k, \delta\Pi_k$ which are independent, the expression in the
braces are zero. Thus we obtain necessary conditions for optimality
as follows.
\begin{gather}
\Pi_{k+1} = F_k^T \Pi_k + h \parenth{M_{k+1}+Bu_{k+1}},\label{eqn:updatePi1}\\
h S(\Pi_k) = F_k J_d - J_dF_k^T,\label{eqn:findf1}\\
R_{k+1} = R_k F_k,\label{eqn:updateR1}\\
u_{k+1} = -B^T\lambda_{k}^2,\label{eqn:ukp}\\
\begin{bmatrix} \lambda_{k}^1 \\ \lambda_{k}^2 \end{bmatrix}
= \begin{bmatrix} \mathcal{A}_{k+1}^T & \mathcal{C}_{k+1}^T \\
\mathcal{B}_{k+1}^T & \mathcal{D}_{k+1}^T\end{bmatrix}
\begin{bmatrix} \lambda_{k+1}^1 \\ \lambda_{k+1}^2
\end{bmatrix},\label{eqn:updatelam}
\end{gather}
where
\begin{align}
\mathcal{A}_k & = F_k^T,\\
\mathcal{B}_k & = hF_k^T\braces{\tr{F_kJ_d}I_{3\times 3}-F_kJ_d}^{-1},\label{eqn:BBk}\\
\mathcal{C}_k & = h\mathcal{M}_{k+1}F_k^T,\\
\mathcal{D}_k & = F_k^T + S(F_k^T\Pi_k)\mathcal{B}_k+h
\mathcal{M}_{k+1}\mathcal{B}_k.
\end{align}
In the above equations, the only implicit part is \refeqn{findf1}.
For a given initial condition $(R_0,\Pi_0,\lambda_0^1,\lambda_0^2)$,
we can find $F_0$ by solving \refeqn{findf1}. Then, $R_1$ is
obtained from \refeqn{updateR1}. Since $u_1=-\lambda_0^2$ by
\refeqn{ukp}, and $M_{1}$ is a function of $R_1$, $\Pi_1$ can be
obtained using \refeqn{updatePi1}. We solve \refeqn{findf1} to
obtain $F_1$ using $\Pi_1$. Finally, $\lambda_1^1,\lambda_1^2$ are
obtained from \refeqn{updatelam}, since
$\mathcal{A}_1,\mathcal{B}_1,\mathcal{C}_1,\mathcal{D}_1$ are
functions of $R_1,\Pi_1,F_1$. This yields a map
$(R_0,\Pi_0,\lambda_0^1,\lambda_0^2)\mapsto(R_1,\Pi_1,\lambda_1^1,\lambda_1^2)$,
and this process can be repeated.

\subsection{Two point boundary value problem}\label{subsec:tpbvp} The necessary
conditions for optimality are given as a 12 dimensional two point
boundary value problem on $T^*\SO$ and its dual space. This problem
is to find
\begin{align*}
\text{Attitude and Angular momentum :} & R_k, \Pi_k,\\
\text{Multiplier variables :} & \lambda^1_k, \lambda^2_k,\\
\text{Control Inputs :} & u_k,
\end{align*}
for $k=\braces{0,\,1,\,\cdots,N}$, to satisfy simultaneously,
\begin{align*}
\text{Equations of motion :} & \refeqn{updatePi1}, \refeqn{findf1}, \refeqn{updateR1},\\
\text{Multiplier equations :} & \refeqn{updatelam},\\
\text{Optimality condition :} & \refeqn{ukp},\\
\text{Boundary conditions :} & R_0,\Pi_0,R_N,\Pi_N.
\end{align*}
An iterative numerical method for the two point boundary value
problem is presented. A nominal solution that satisfies some of the
above conditions is chosen, and this nominal solution is updated by
successive linearization so that the remaining conditions are also
satisfied as the process converges.

We use a neighboring extremal method~\cite{bk:bryson}. A nominal
solution satisfies all of the necessary conditions except the
boundary conditions. The neighboring extremal method is
characterized as an iterative algorithm for improving estimates of
the unspecified multiplier initial conditions so as to satisfy the
specified terminal boundary conditions in the limit. This is
sometimes referred to as a shooting method. The optimality condition
\refeqn{ukp} is substituted into the equations of motion and the
multiplier equations. The sensitivities of the specified terminal
boundary conditions with respect to the unspecified initial
multiplier conditions can be calculated by direct numerical
differentiation, or they can be obtained by a linear analysis. The
main advantage of the neighboring extremal method is that the number
of iteration variables is small. It is equal to the dimension of the
equations of motion. The difficulty is that the extremal solutions
are sensitive to small changes in the unspecified initial multiplier
values. Therefore, it is important to compute the sensitivities
accurately.

We use linear analysis to compute the sensitivities. The sensitivity
model is defined at the Lie algebra level as presented
in~\cite{pro:acc06}. It is natural to define the sensitivity model
in the Lie algebra, since the Lie algebra is a linear vector space.
The resulting sensitivity model is global, and it has the same
dimension as the Lie group. The sensitivity derivatives in the Lie
algebra are related to the original Lie group by the exponential
map.

Using the perturbation models defined in \refeqn{delPik},
\refeqn{delRk}, the linearized equations of motion for the attitude
dynamics can be written as
\begin{gather}
\begin{bmatrix} \zeta_{k+1} \\ \delta\Pi_{k+1} \end{bmatrix}
= \begin{bmatrix} \mathcal{A}_{k} & \mathcal{B}_{k} \\
\mathcal{C}_{k} & \mathcal{D}_{k}\end{bmatrix}
\begin{bmatrix} \zeta_{k} \\ \delta\Pi_{k} \end{bmatrix}
-\begin{bmatrix} 0_{3\times 3}\\hBB^T\end{bmatrix} \delta
\lambda_{k}^2 .\label{eqn:delx}
\end{gather}
Note that the homogeneous part of \refeqn{delx} is equivalent to
equations that are the dual of \refeqn{updatelam}. The variation of
the equations of motion is equivalent to the dual of the multiplier
equations, so the variation of the multiplier equations is
equivalent to the second variation of the attitude dynamics
equations. The linearized equations of motion for the multipliers
can be obtained as
\begin{gather}
\begin{bmatrix} \delta\lambda_{k+1}^1 \\ \delta\lambda_{k+1}^2 \end{bmatrix}
=-\parenth{A^{11}_{k+1}}^{-T}A^{21}_{k+1}A^{11}_{k}
\begin{bmatrix} \zeta_{k} \\ \delta\Pi_{k} \end{bmatrix}
+\parenth{A^{11}_{k+1}}^{-T}\parenth{I_{3\times
3}-A^{21}_{k+1}A^{12}_{k}}
\begin{bmatrix} \delta\lambda^1_{k} \\ \delta\lambda^2_{k} \end{bmatrix}
,\label{eqn:dellambda}
\end{gather}
where $A_k^{ij}\in\Re^{6\times 6},i,j\in\braces{1,2}$ are defined in
the Appendix.

Using \refeqn{delx} and \refeqn{dellambda}, the sensitivity
derivatives of the attitude, angular momentum, and the multipliers
can be written as
\begin{align}
\begin{bmatrix}x_{k+1}\\\delta\lambda_{k+1}\end{bmatrix}
= A_k
\begin{bmatrix}x_{k}\\\delta\lambda_{k}\end{bmatrix},\label{eqn:sen}
\end{align}
where $x_k=[\zeta_k^T,\delta\Pi_k^T]^T$,
$\delta\lambda_k=[\delta\lambda_k^{1,T},\delta\lambda_k^{2,T}]^T\in\Re^6$,
and $A_k\in\Re^{12\times 12}$. The solution of \refeqn{sen} is given
by
\begin{align}
\begin{bmatrix}x_{N}\\\delta\lambda_N\end{bmatrix}& =
\left(\prod_{k=0}^{N-1}A_k\right)
\begin{bmatrix}x_{0}\\\delta\lambda_0\end{bmatrix},\nonumber\\
& \triangleq
\begin{bmatrix} \Phi_{11} & \Phi_{12} \\ \Phi_{21} & \Phi_{22}\end{bmatrix}
\begin{bmatrix}x_{0}\\\delta\lambda_0\end{bmatrix}.
\end{align}
For the given two point boundary value problem $x_0=0$ since the
initial attitude and the initial angular momentum are given, and
$\lambda_N$ is free.  Then, we obtain
\begin{align}
x_N = \Phi_{12} \delta\lambda_0
\end{align}
The unspecified initial multipliers are $\lambda_0$, and the
specified terminal boundary conditions are the terminal attitude
$R_N^d$ and the terminal angular momentum $\Pi_N^d$. Thus
$\Phi_{12}$ represents the sensitivity of the specified terminal
boundary conditions with respect to the unspecified initial
multipliers. Using this sensitivity, an initial guess of the
unspecified initial conditions is iterated to satisfy the specified
terminal conditions in the limit.

Any type of Newton iteration can be applied to this problem using
the sensitivity derivative as a gradient. The procedure uses a line
search with backtracking algorithm, referred to as Newton-Armijo
iteration in~\cite{bk:kelley}. The procedure is summarized in Table
\ref{tab:iter}, where $i$ is the iteration index, and
$\epsilon_S,\alpha\in\Re$ are a stopping criterion and a scaling
factor, respectively. The outer loop finds a search direction by
computing the sensitivity derivatives, and the inner loop performs a
line search along the obtained direction; the error in satisfaction
of the terminal boundary condition is determined on each iteration.

\section{Numerical Computations}\label{sec:sim}
Numerical results are given for two optimal attitude control
problems; optimal attitude control of an underactuated 3D pendulum
and optimal attitude control of a fully actuated spacecraft on a
circular orbit.

\subsection{3D Pendulum}
A 3D pendulum is a rigid body supported at a frictionless pivot
acting under the influence of uniform gravity~\cite{pro:nalin}. The
gravity potential, acting in the vertical or $e_3$ direction, is
given by
\begin{align}
U(R) & = - mg e_3^T R \rho,\label{eqn:pendU}
\end{align}
where $m\in\Re$ is the mass of the pendulum, $g\in\Re$ is the
gravitational acceleration, and $\rho\in\Re^3$ represents a vector
from the pivot point to the mass center of the pendulum in the body
fixed frame. The pendulum model is shown in \reffig{body}.(a) with
the pivot located at the origin, and we assume that the pendulum is
axially symmetric. The gravity moment and its variations are given
by
\begin{align*}
M & = mg \rho \times R^T e_3,\\
\delta M & = \mathcal{M}\zeta = mg S(\rho) S(R^T e_3) \zeta.
\end{align*}
There are two equilibrium manifolds; a hanging equilibrium manifold
when $R\rho=\norm{\rho}e_3$, an inverted equilibrium manifold when
$R\rho=-\norm{\rho}e_3$.

The properties of the axially symmetric pendulum are given by
$J=\rm{diag}\bracket{0.156,\, 0.156,\, 0.3}\mathrm{kg\,m^2}$, $m=1
\mathrm{kg}$, and $\rho=\bracket{0, 0, \frac{3}{4}} \mathrm{m}$. We
assume that the component of the control input along the axis of
symmetry is zero; this corresponds to an underactuated 3D pendulum.
The corresponding input matrix is given by
\begin{align*}
B=\begin{bmatrix} 1 & 0 \\ 0 & 1\\ 0 & 0\end{bmatrix}.
\end{align*}

Two types of boundary conditions are considered. The first maneuver
is to transfer the 3D pendulum from a hanging equilibrium to an
inverted equilibrium. The second maneuver is a
$180\,\mathrm{degree}$ rotation about the uncontrolled axis of
symmetry starting in a hanging equilibrium. The terminal attitude
also lies in the hanging equilibrium manifold. Each maneuver is
completed in $1\,\mathrm{sec}$. The time step size is
$h=0.001\,\mathrm{sec}$ and the number of integration steps is
$N=1000$. The corresponding boundary conditions are given by
\renewcommand{\theenumi}{\roman{enumi}}
\renewcommand{\labelenumi}{(\theenumi)}
\begin{enumerate}
\item Rotation from a hanging equilibrium to an inverted
equilibrium.
\begin{gather*}
R_{0}=I_{3\times 3},\quad R_{N}^d=\begin{bmatrix}
     0&     1&     0\\
     1&     0&     0\\
     0&     0&     -1\end{bmatrix},\\
\Pi_0=0_{3\times 1},\quad\Pi_N^d=0_{3\times 1}.
\end{gather*}
\item Rotation from one hanging equilibrium to another hanging equilibrium.
\begin{gather*}
R_{0}=I_{3\times 3},\quad R_{N}^d=\mathrm{diag}[-1,-1,1],\\
\Pi_0=0_{3\times 1},\quad\Pi_N^d=0_{3\times 1}.
\end{gather*}

\end{enumerate}

\begin{figure*}[h]
    \centerline{\subfigure[Attitude maneuver][]{
    \includegraphics[width=0.45\textwidth]{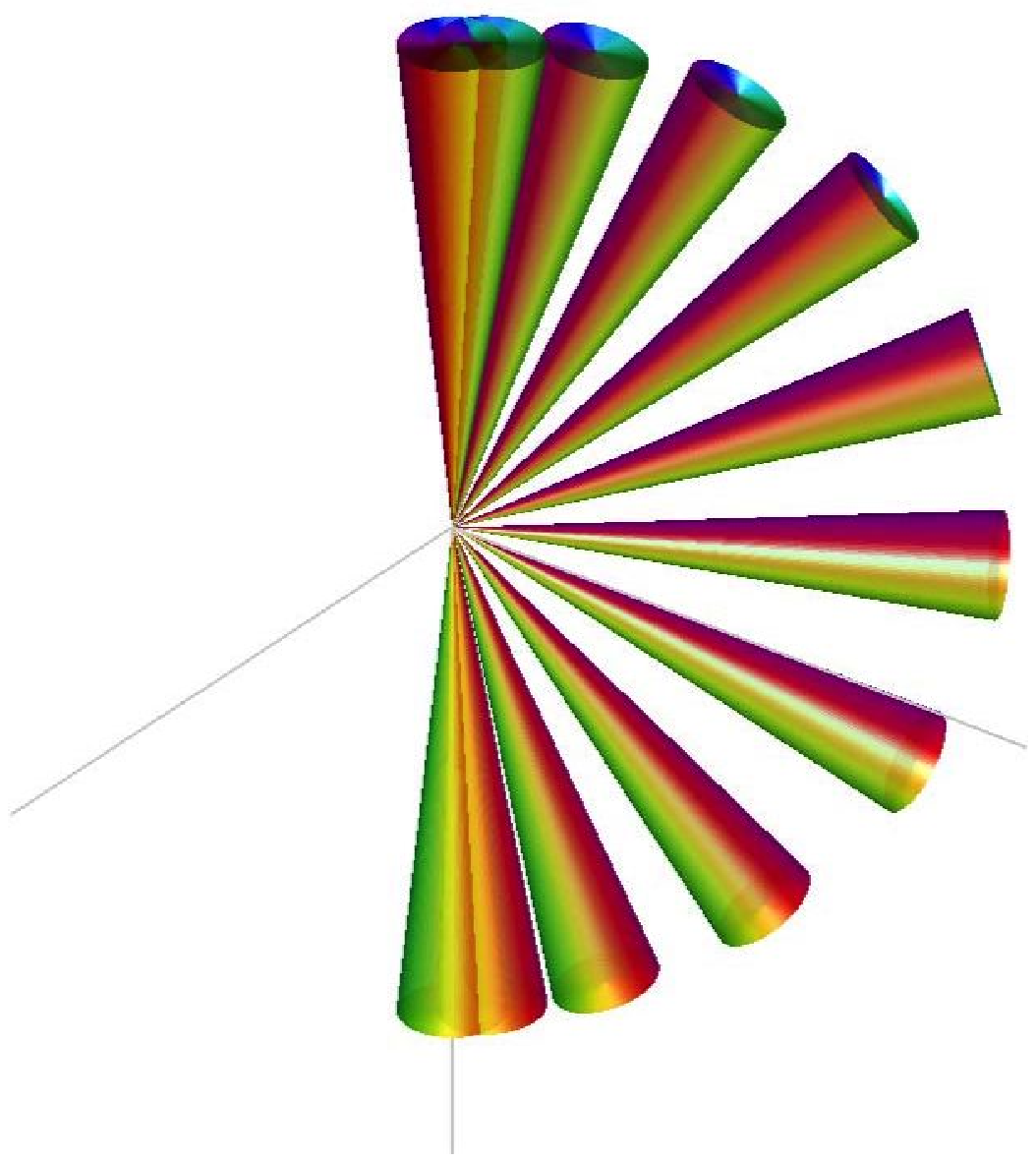}}
    \hfill
    \subfigure[Control input $u$][]{
    \includegraphics[width=0.45\textwidth]{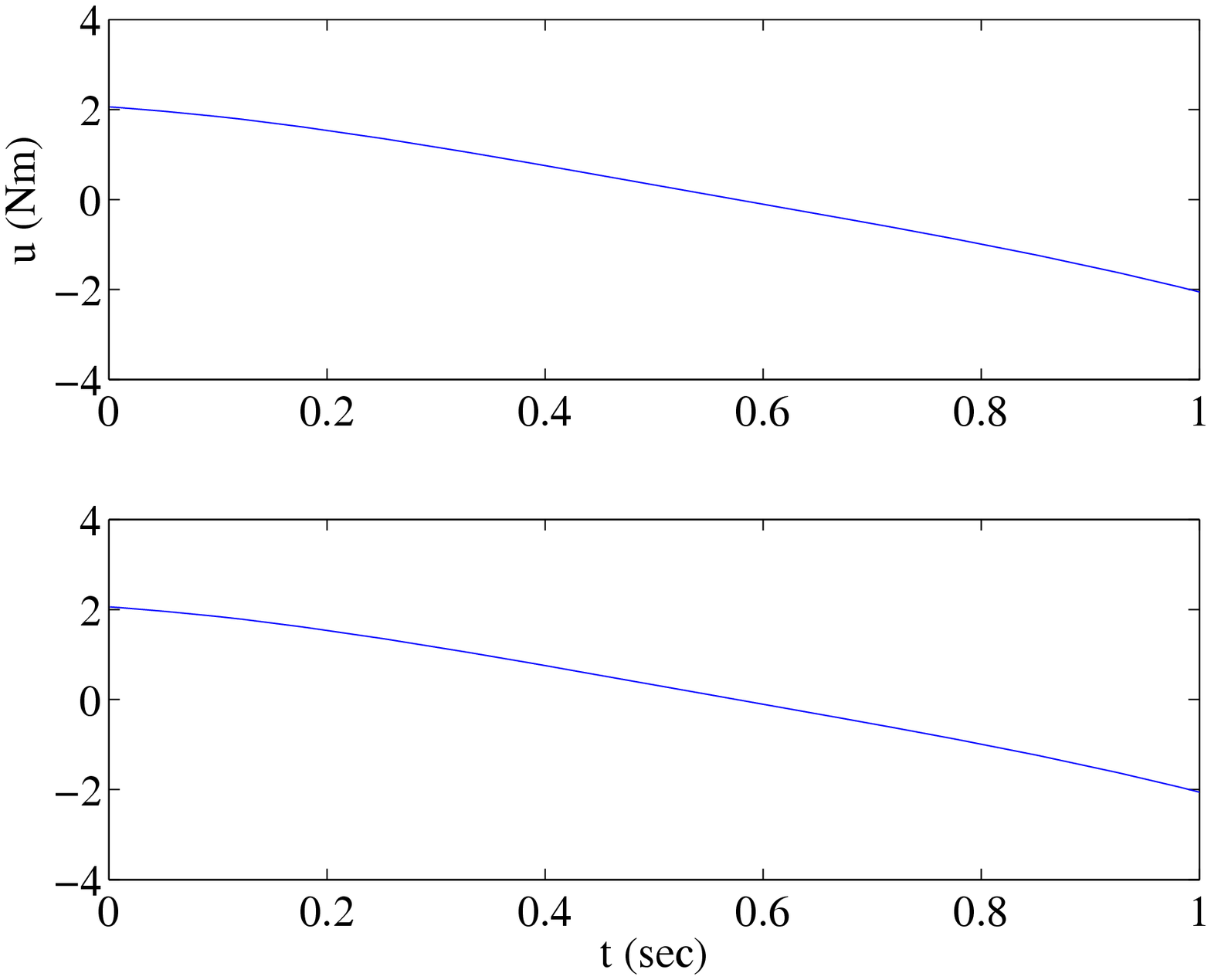}}
    }
    \vspace*{1cm}
    \centerline{\subfigure[Angular velocity $\Omega$][]{
    \includegraphics[width=0.45\textwidth]{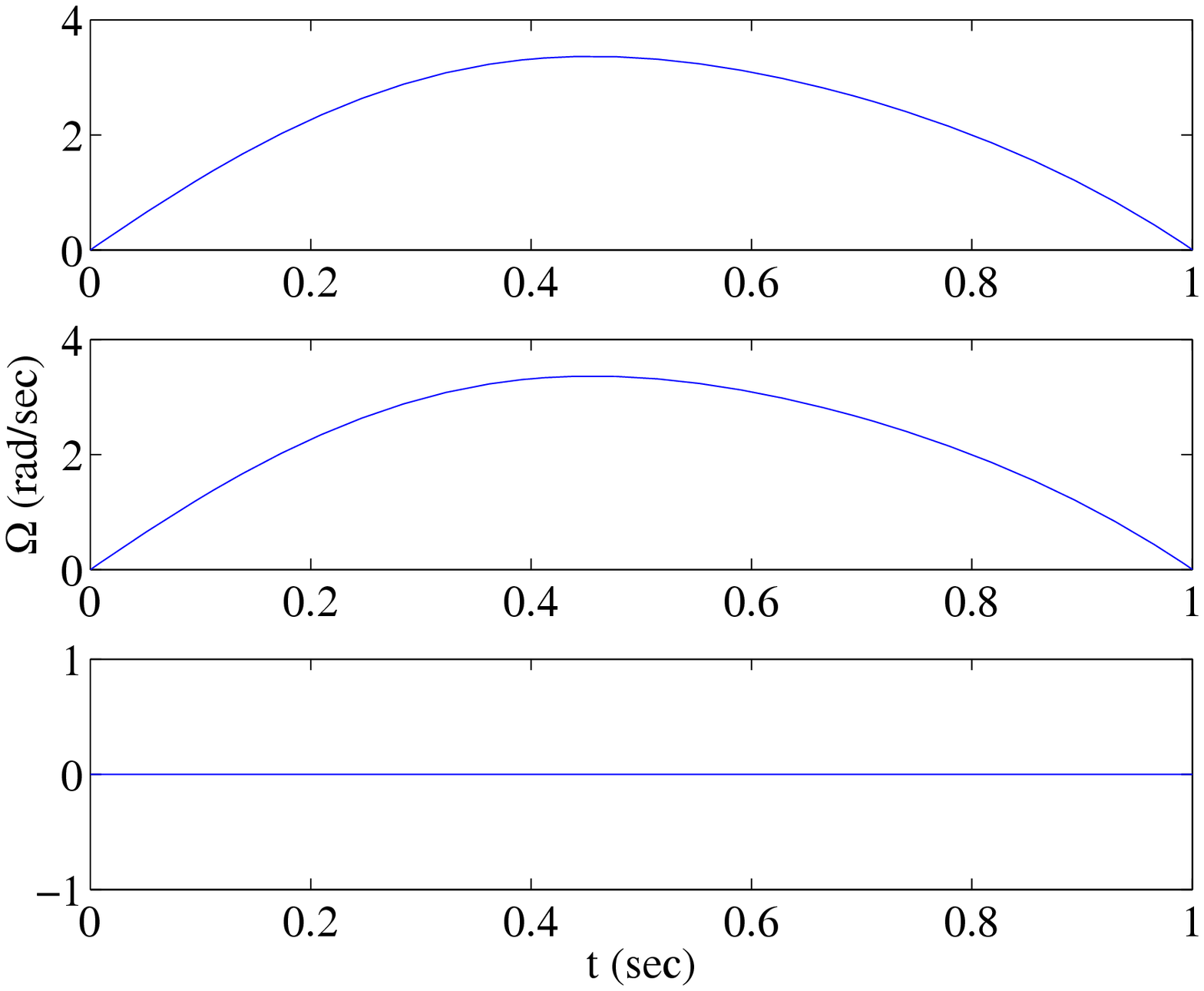}}
    \hfill
    \subfigure[Convergence rate][]{
    \includegraphics[width=0.45\textwidth]{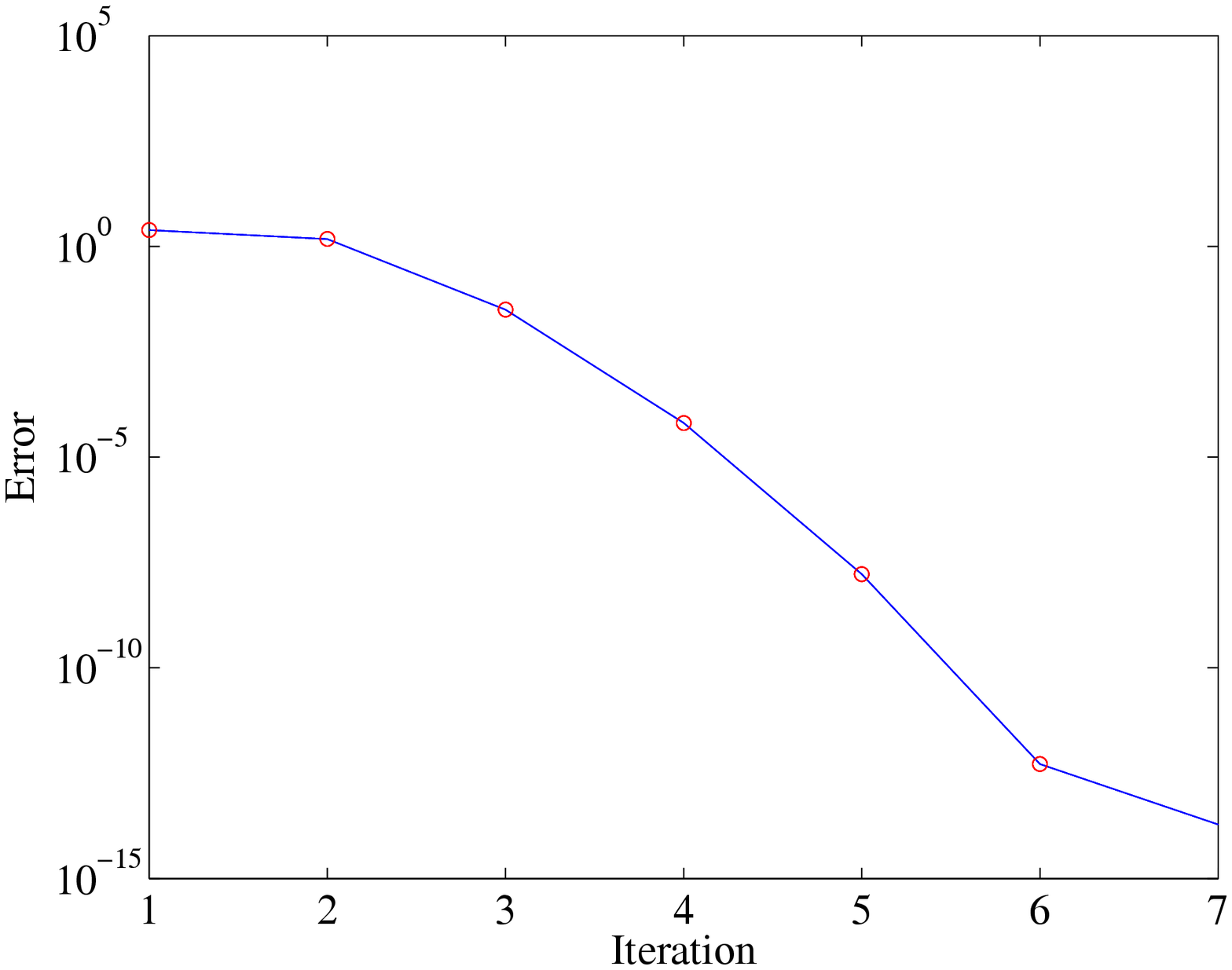}}
    }
    \caption[Case (i): 3D pendulum rotation from a hanging equilibrium to an inverted
equilibrium]{}\label{fig:pendi}
\end{figure*}

\begin{figure*}[h]
    \centerline{\subfigure[Attitude maneuver][]{
    \includegraphics[width=0.45\textwidth]{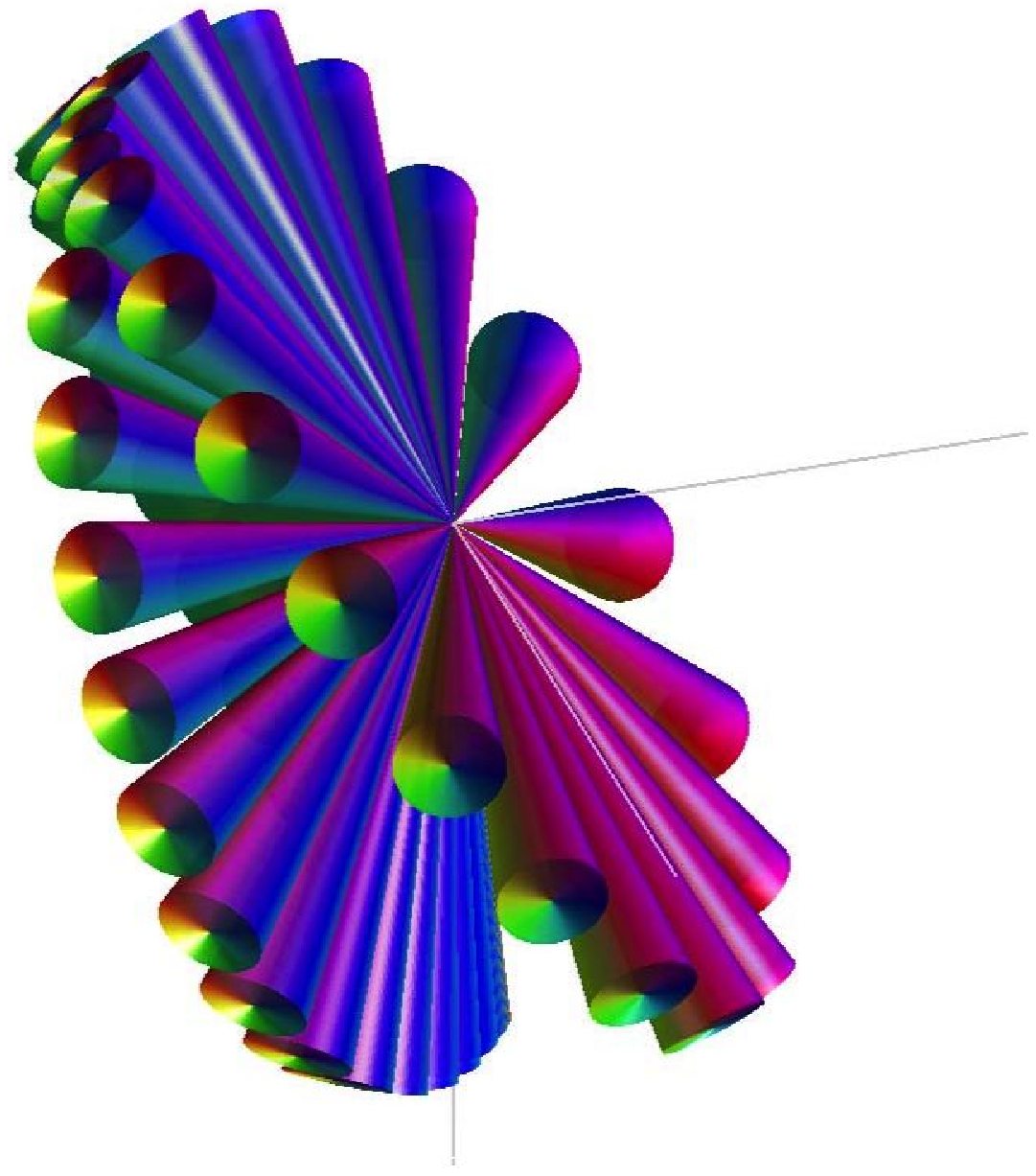}}
    \hfill
    \subfigure[Control input $u$][]{
    \includegraphics[width=0.45\textwidth]{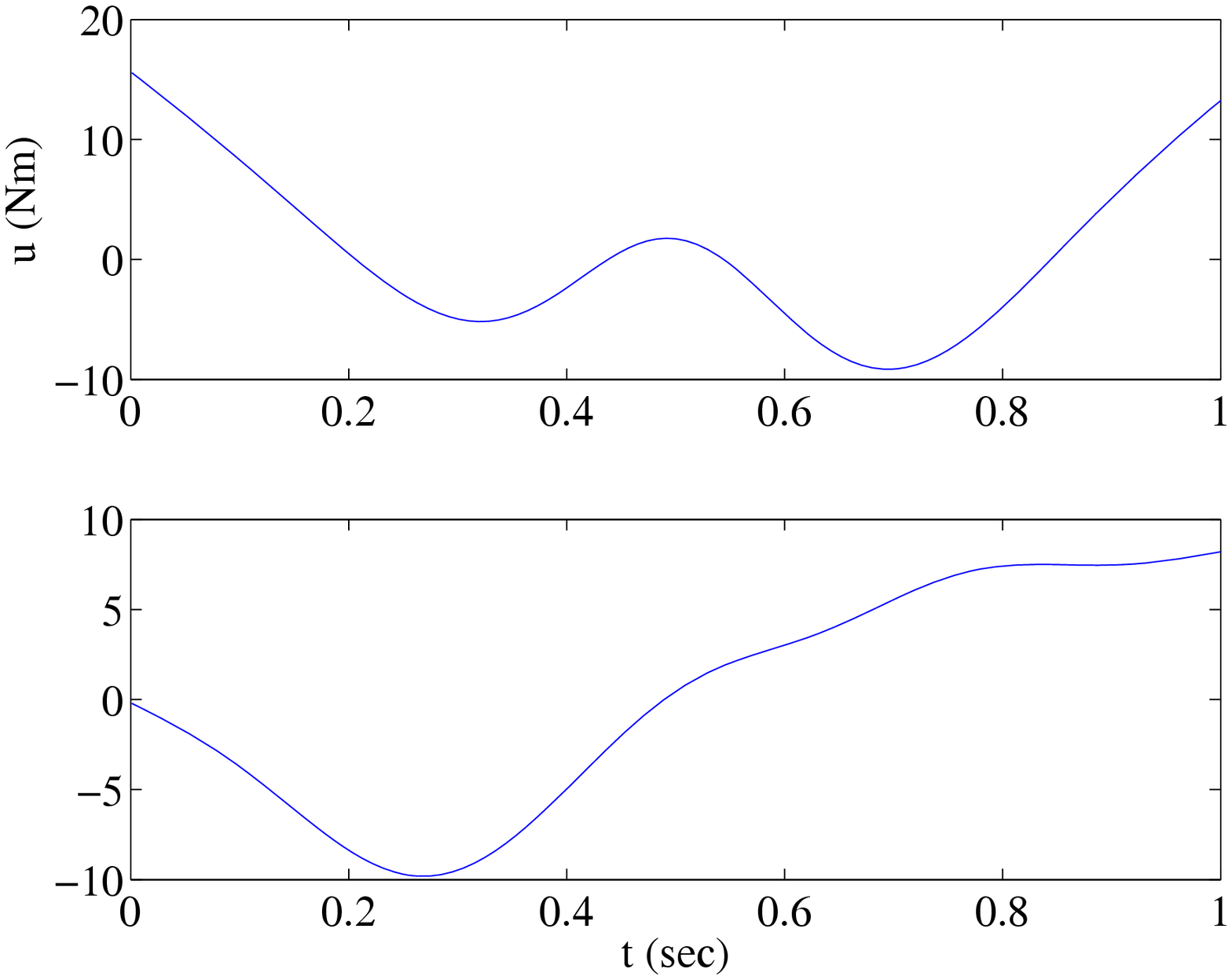}}
    }
    \vspace*{1cm}
    \centerline{\subfigure[Angular velocity $\Omega$][]{
    \includegraphics[width=0.45\textwidth]{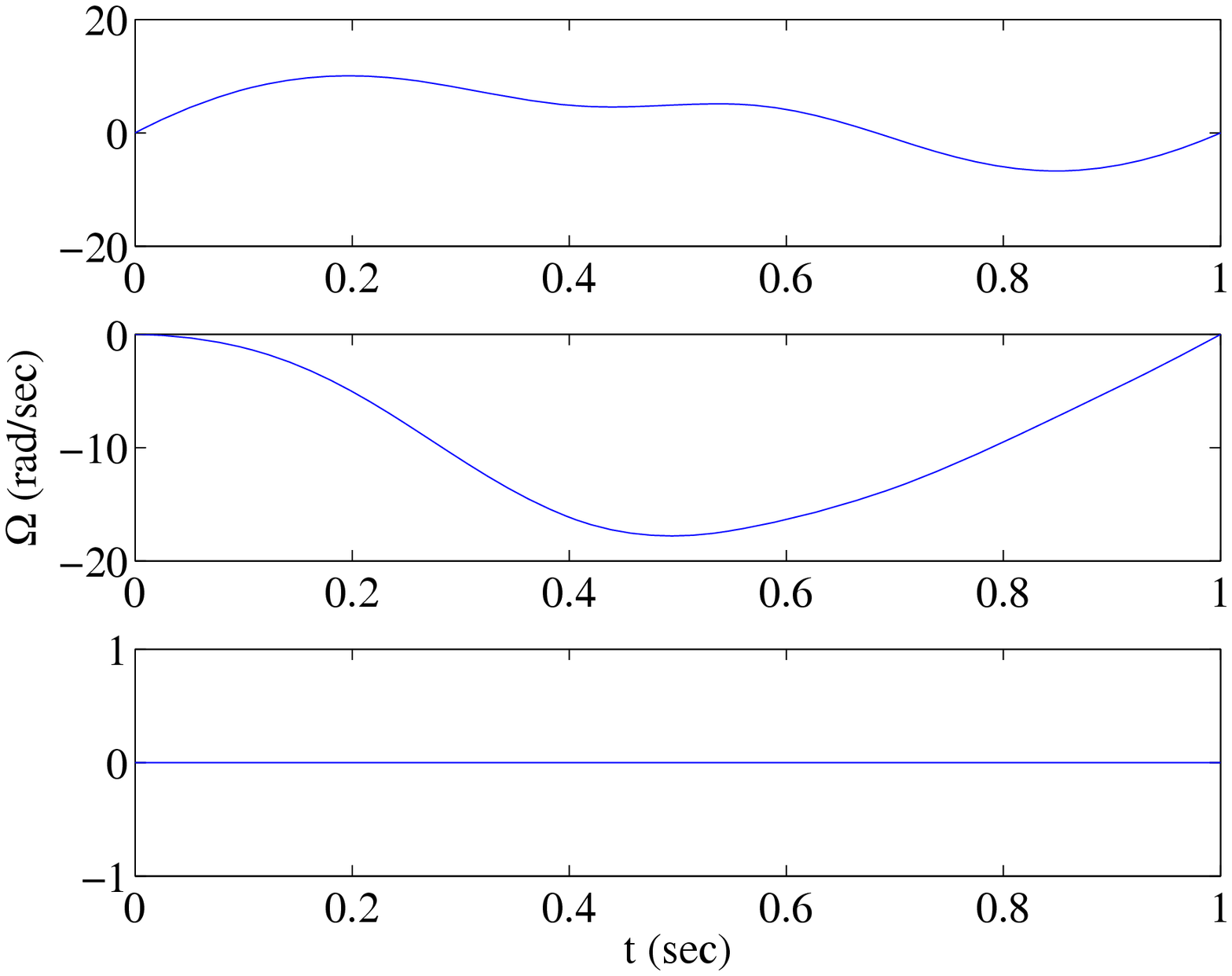}}
    \hfill
    \subfigure[Convergence rate][]{
    \includegraphics[width=0.45\textwidth]{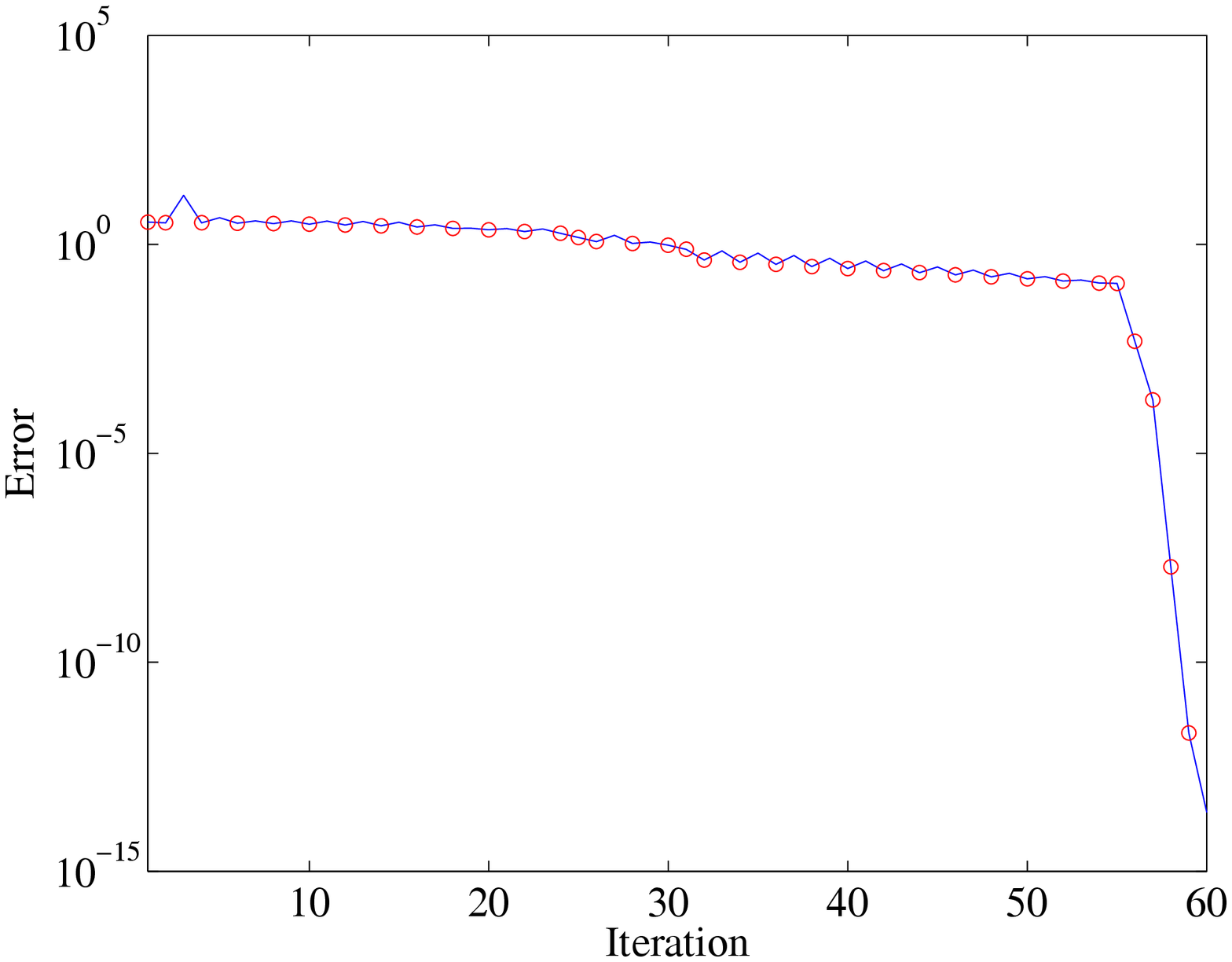}}
    }
    \caption[Case (ii): 3D pendulum rotation from one hanging equilibrium to another hanging equilibrium]{}\label{fig:pendii}
\end{figure*}

The optimized performance index and the violation of the constraints
are given in Table \ref{tab:res}. The terminal boundary conditions
are satisfied at the level of machine precision for both cases.
Figures \ref{fig:pendi} and \ref{fig:pendii} show snapshots of the
attitude maneuvers, the control input history, and the angular
velocity response. (Simple animations which show these optimal
attitude maneuvers of the 3D pendulum can be found at
\url{http://www.umich.edu/~tylee}.) The third component of the
angular velocity is constant; this is a conservation property of the
controlled axially symmetric 3D pendulum.

The optimized attitude maneuver of the first case is an eigen-axis
rotation about the fixed axis;
$[\frac{\sqrt{2}}{2},\,\frac{\sqrt{2}}{2},\,0]$. In the second case,
the rotation about the axis of symmetry is induced from control
moments about the first and second body fixed axes. The resulting
attitude maneuver is more complicated, and it requires larger
control inputs.

Figures \ref{fig:pendi}.(d) and \ref{fig:pendii}.(d) show the
violation of the terminal boundary condition according to the number
of iterations in a logarithm scale. The circles denote outer
iterations to compute the sensitivity derivatives. For all cases,
the initial guesses of the unspecified initial multiplier are
arbitrarily chosen such that the initial trial of control inputs is
close to zero throughout the maneuver time. The error in
satisfaction of the terminal boundary condition of the first case
converges quickly to machine precision; only 7 iterations are
required. A longer number of  iterations is required in the second
case, but the error converges exponentially to machine precision
after the solution is close to the local minimum at the 55th
iteration. These convergence rates show the quadratic convergence
property of Newton iteration. The convergence rates are dependent on
the numerical accuracy of the sensitivity derivatives.

\subsection{Spacecraft on a Circular Orbit}
We consider a spacecraft on a circular orbit about a large central
body, including gravity gradient effects~\cite{bk:wie}. The
spacecraft model is shown at \reffig{body}.(b). The attitude of the
spacecraft is represented with respect to the local vertical local
horizontal (LVLH) axes. The gravity potential is given by
\begin{align*}
U(R)=-\frac{GM}{r_0}-\frac{1}{2}\omega_0^2\parenth{\tr{J}-3e_3^TRJR^Te_3},
\end{align*}
where $G\in\Re$ is the gravitational constant, $M\in\Re$ is the mass
of the central body, $r_0\in\Re$ is the orbital radius, and
$\omega_0=\sqrt{\frac{GM}{r_0^3}}\in\Re$ is the orbital angular
velocity. The gravity moment and its variations are given by
\begin{align*}
M & = 3\omega_0^2 R^Te_3\times JR^Te_3,\\
\delta M & = \mathcal{M}\zeta =3\omega_0^2 \Big[
-S(JR^Te_3)S(R^Te_3)+ S(R^Te_3) J S(R^Te_3)\Big]\zeta
\end{align*}

There are 24 distinct relative equilibria for which the
principal axes are exactly aligned with the LVLH axes, and the
spacecraft angular velocity is identical to the orbital angular
velocity of the LVLH coordinate frame.

We assume that the spacecraft is fully actuated. The corresponding
input matrix is $B=I_{3\times 3}$. The mass, length and time
dimensions are normalized by the mass of the spacecraft, a size
scale factor of the spacecraft, and the orbital angular velocity
$\omega_0\in\Re$, respectively. The mass property of the spacecraft
is chosen as $J=\mathrm{diag}\bracket{1,\,2.8,\,2}$.

\begin{figure*}[h]
    \centerline{\subfigure[Attitude maneuver][]{
    \includegraphics[width=0.45\textwidth]{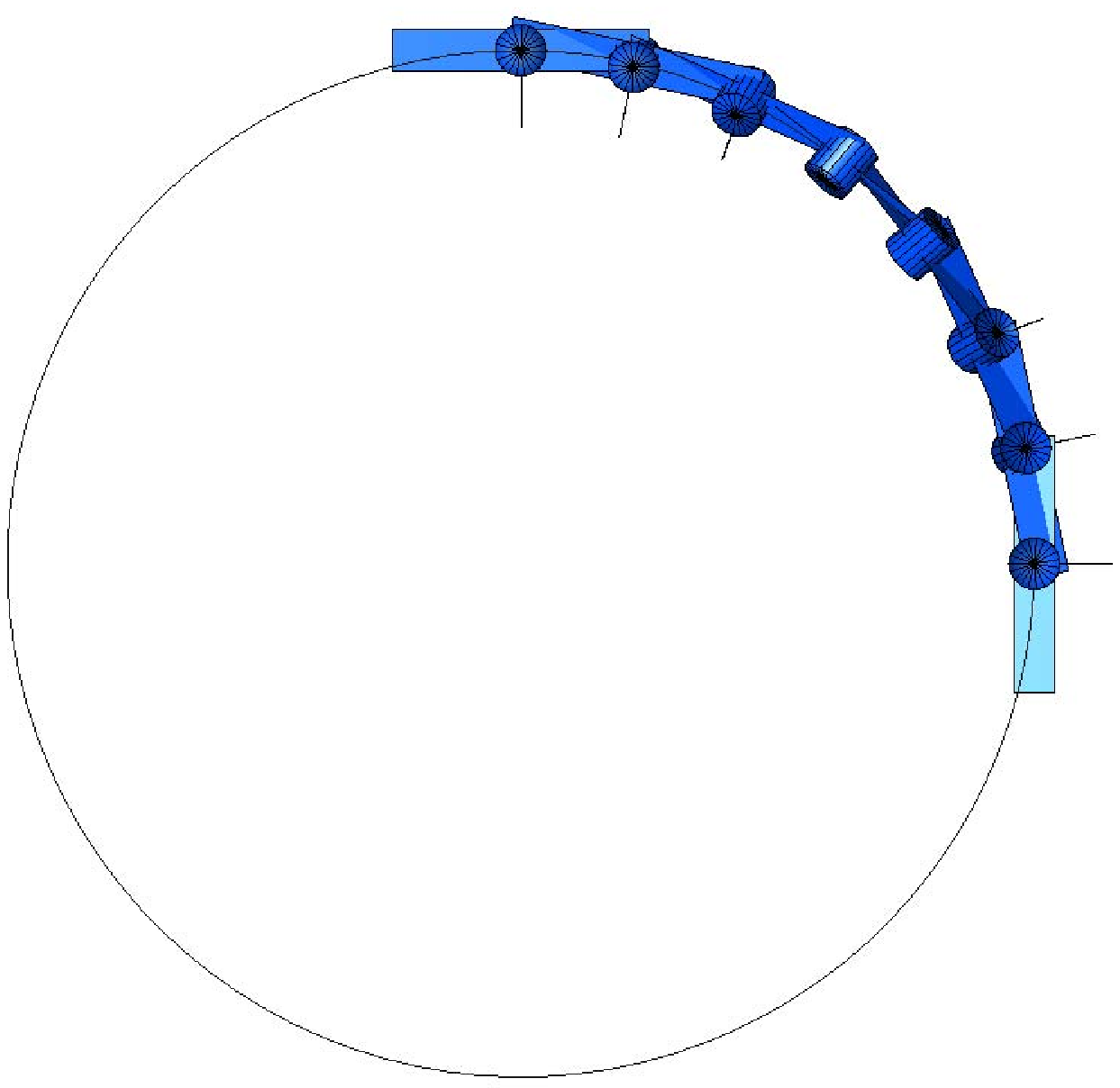}}
    \hfill
    \subfigure[Control input $u$][]{
    \includegraphics[width=0.45\textwidth]{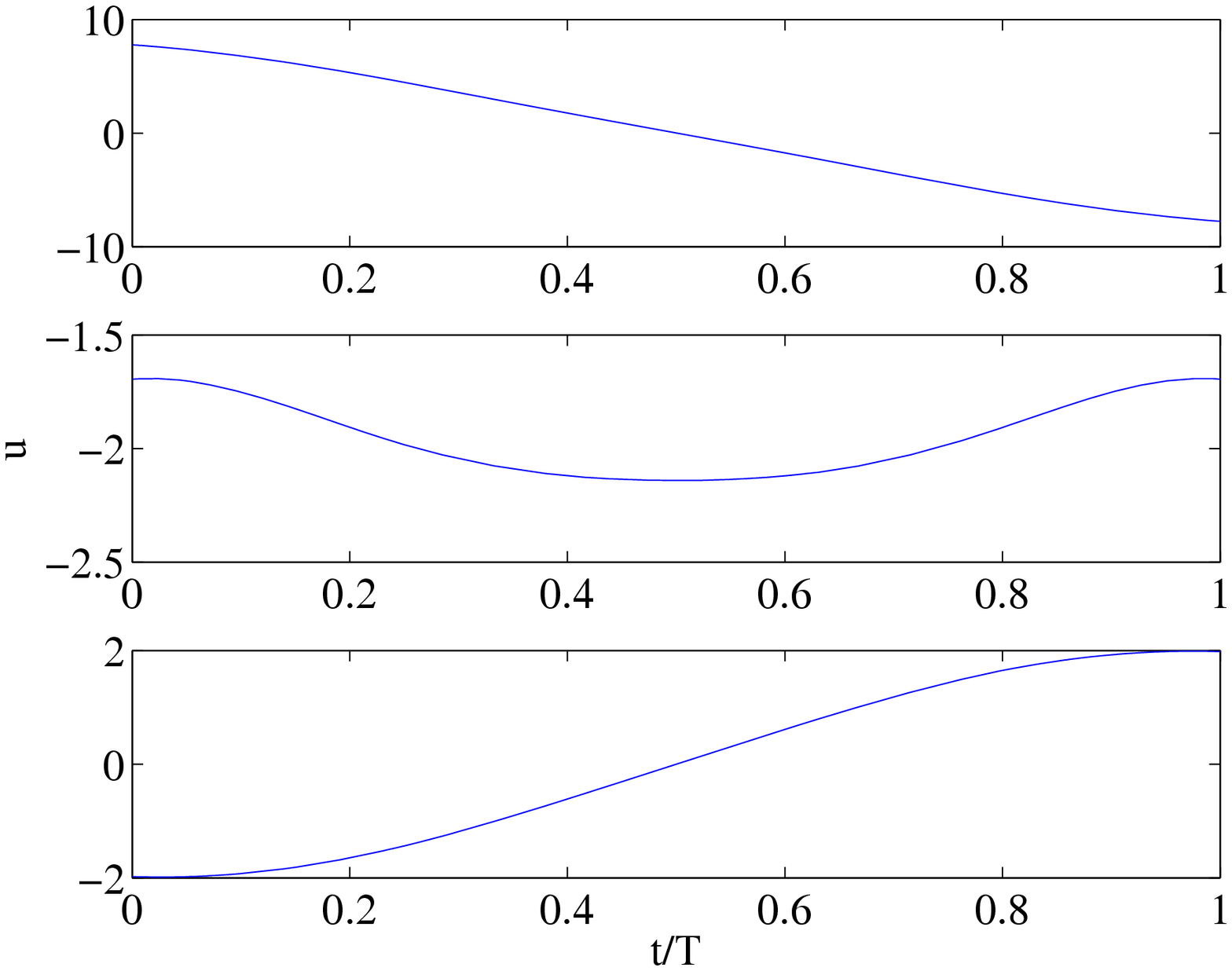}}
    }
    \vspace*{1cm}
    \centerline{\subfigure[Angular velocity $\Omega$][]{
    \includegraphics[width=0.45\textwidth]{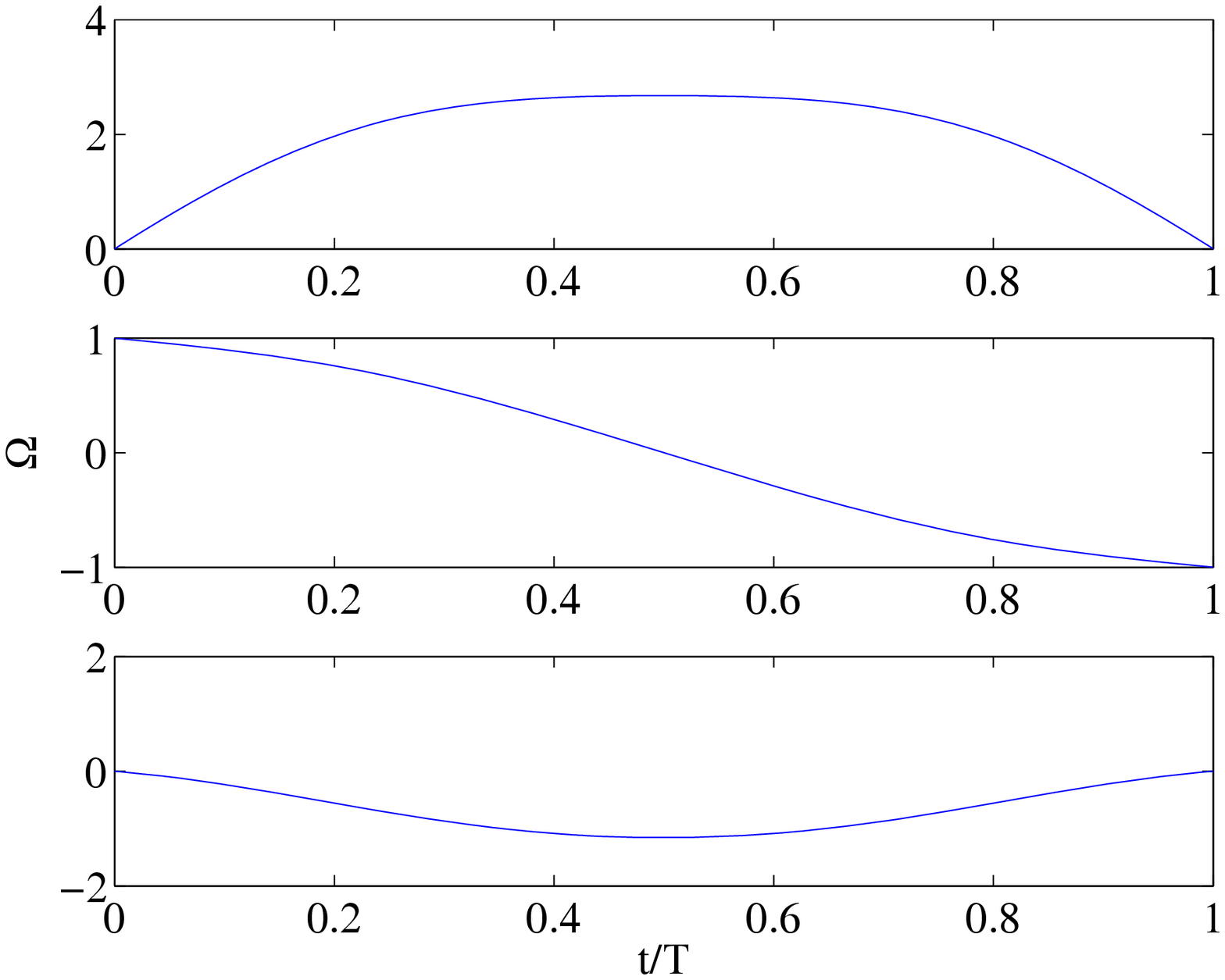}}
    \hfill
    \subfigure[Convergence rate][]{
    \includegraphics[width=0.45\textwidth]{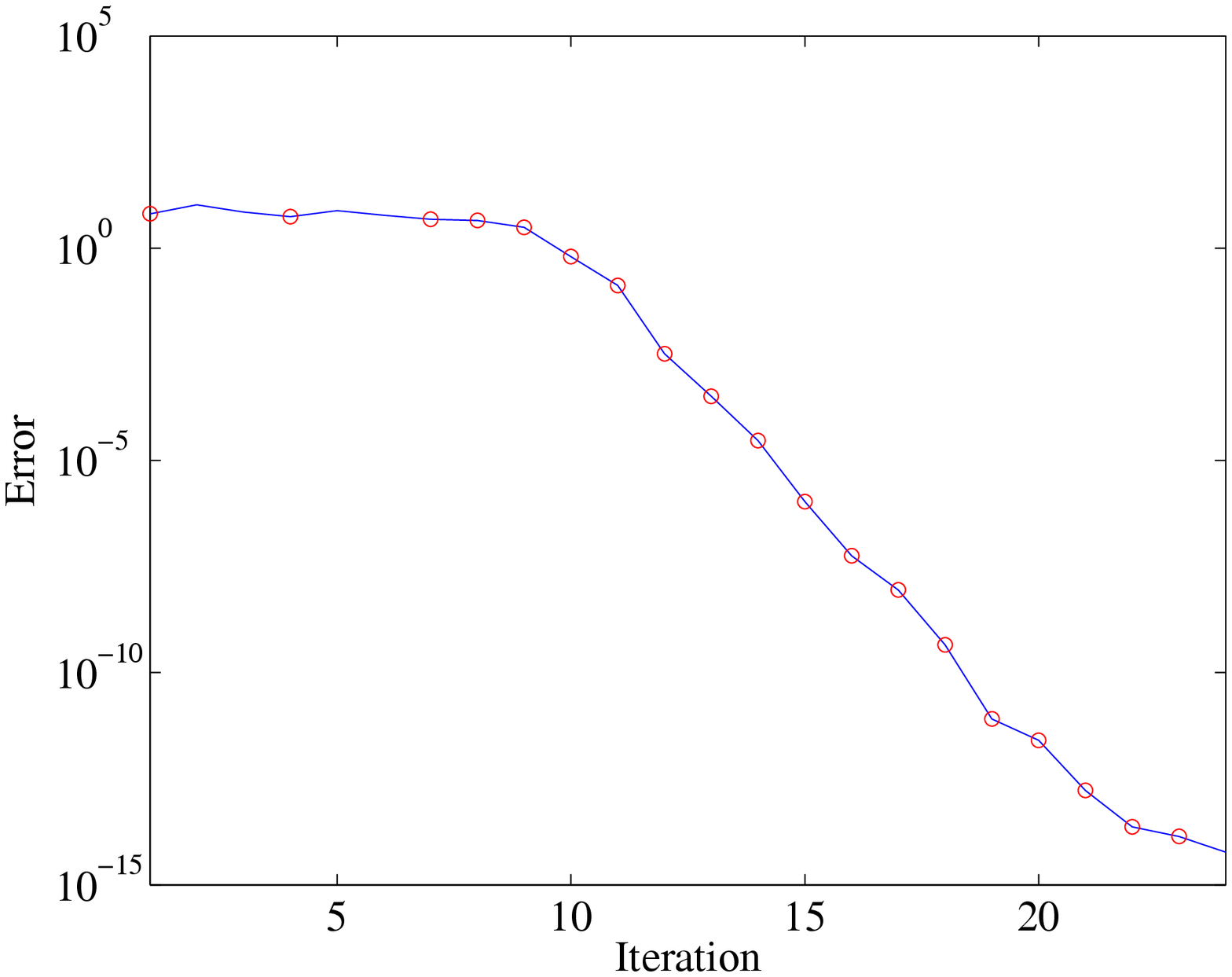}}
    }
    \caption[Case (iii): Spacecraft rotation maneuver about the LVLH tangential axis]{}\label{fig:sati}
\end{figure*}

\begin{figure*}[h]
    \centerline{\subfigure[Attitude maneuver][]{
    \includegraphics[width=0.45\textwidth]{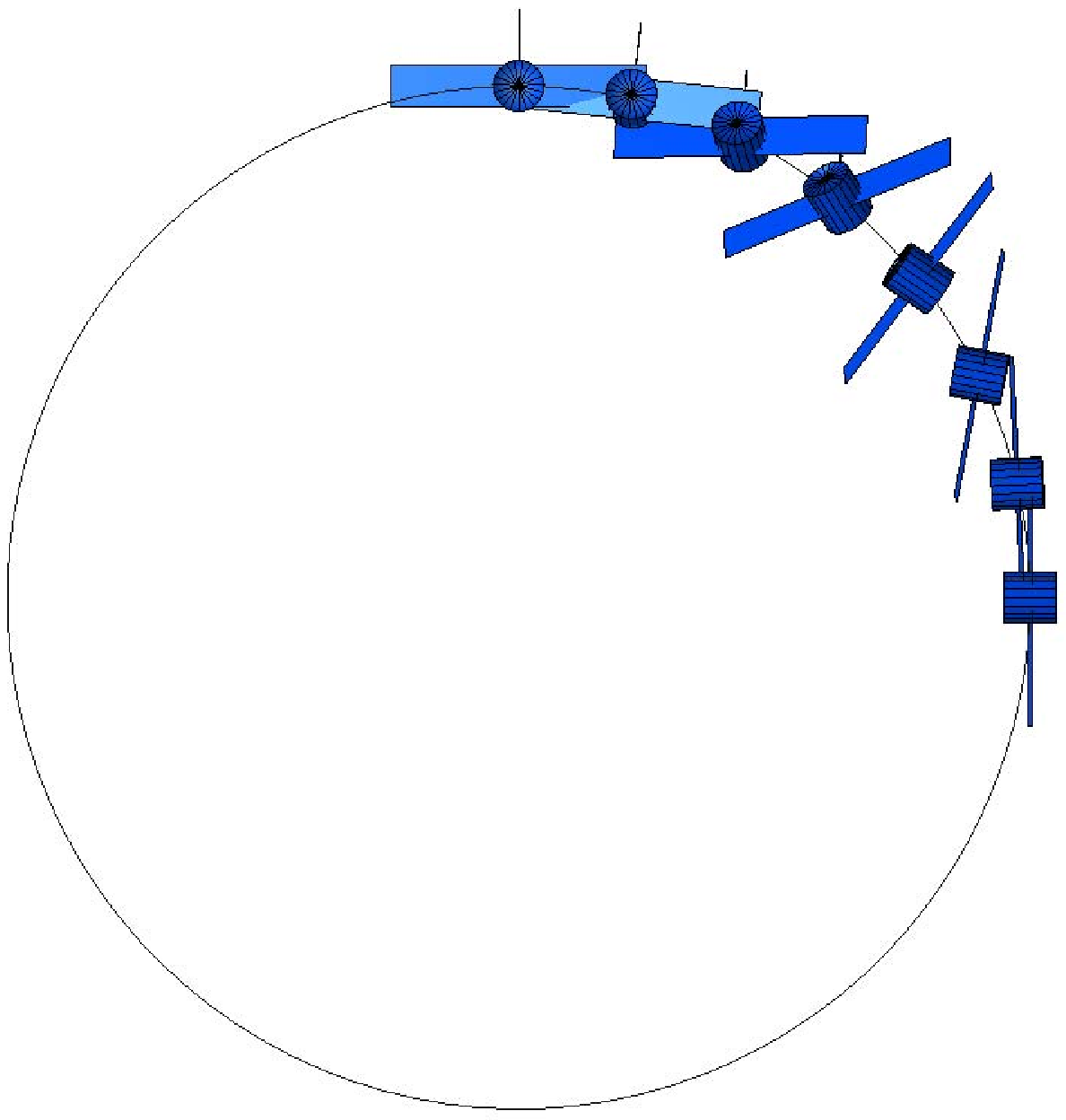}}
    \hfill
    \subfigure[Control input $u$][]{
    \includegraphics[width=0.45\textwidth]{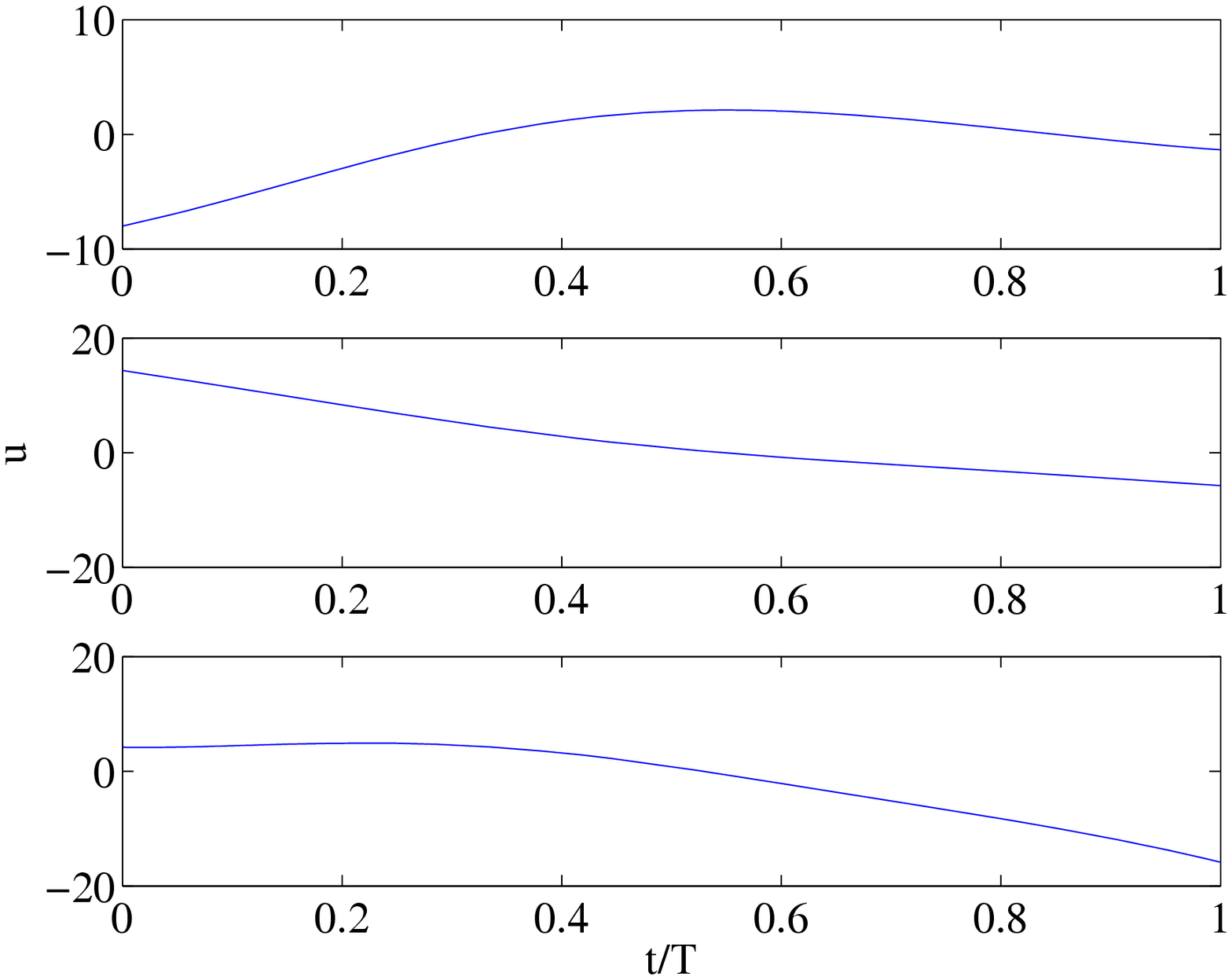}}
    }
    \vspace*{1cm}
    \centerline{\subfigure[Angular velocity $\Omega$][]{
    \includegraphics[width=0.45\textwidth]{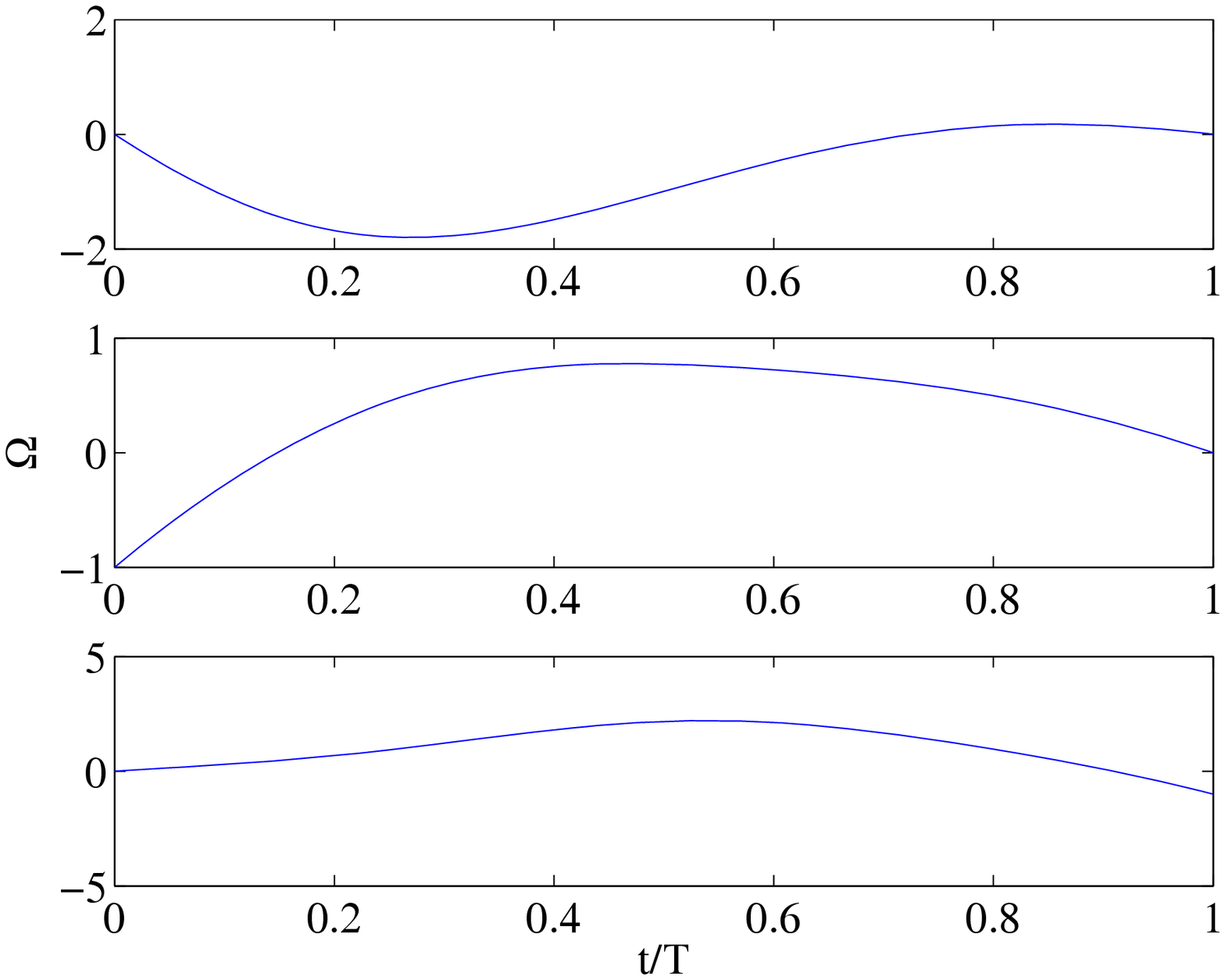}}
    \hfill
    \subfigure[Convergence rate][]{
    \includegraphics[width=0.45\textwidth]{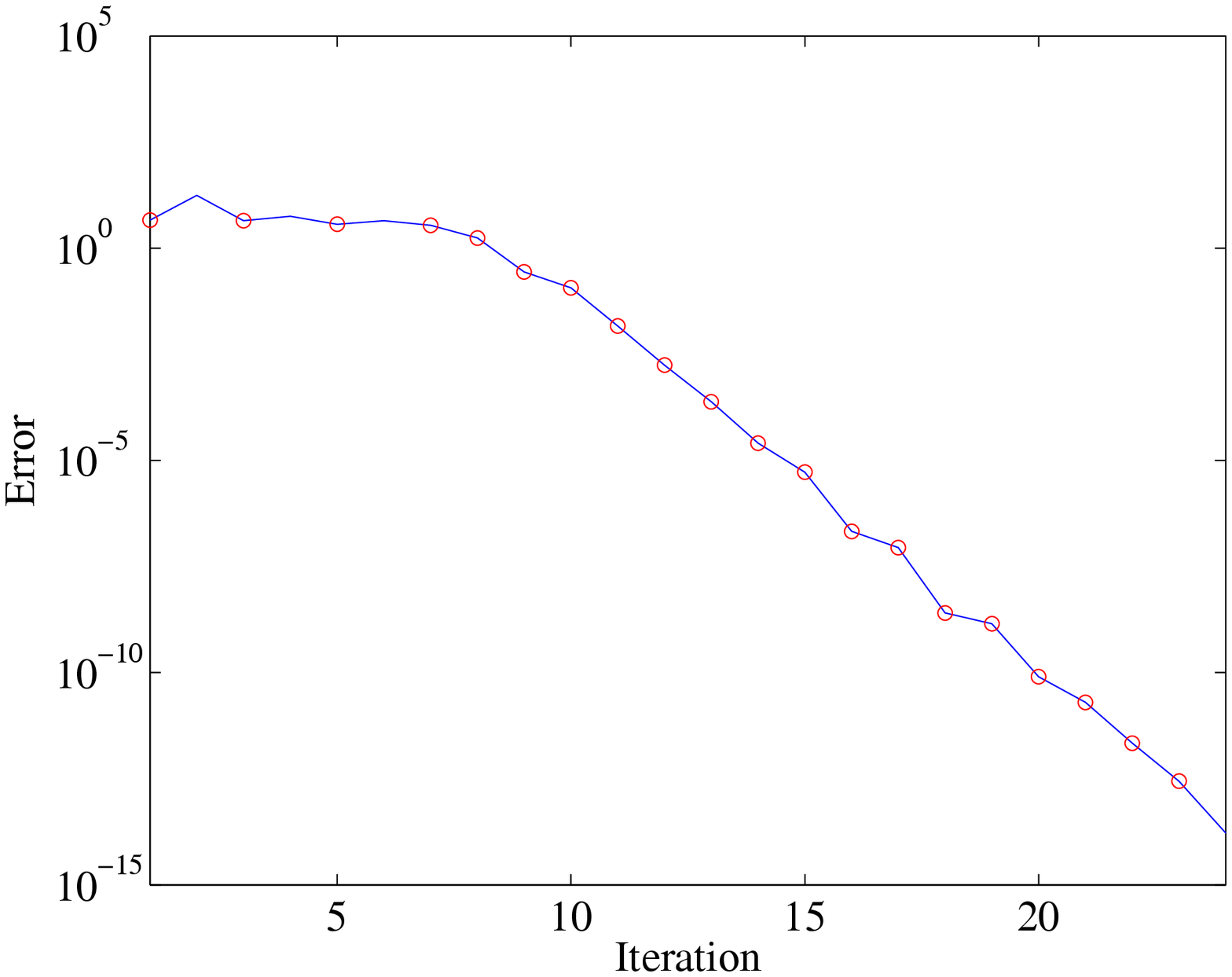}}
    }
    \caption[Case (iv): Spacecraft rotation maneuver about the LVLH tangential and normal axes]{}\label{fig:satii}
\end{figure*}

Two boundary conditions are considered. Each maneuver is a large
attitude change completed in a quarter of the orbit,
$T_f=\frac{\pi}{2}$. The time step size is $h=0.001$ and the number
of integration step is $N=1571$. The terminal angular momentum is
chosen such that the terminal attitude is maintained after the
maneuver.
\renewcommand{\theenumi}{\roman{enumi}}
\renewcommand{\labelenumi}{(\theenumi)}
\begin{enumerate}
\setcounter{enumi}{2}
\item Rotation maneuver about the LVLH axis $e_1$:
\begin{gather*}
R_{0}=I_{3\times 3},\quad
R_{N}^d=\mathrm{diag}\bracket{1, -1, -1},\\
\Pi_0=\omega_0 J R^{T}_0e_2,\quad\Pi_N^d=\omega_0JR^{d,T}_Ne_2.
\end{gather*}
\item Rotation maneuver about the LVLH axes $e_1$ and $e_2$:
\begin{gather*}
R_{0}=\mathrm{diag}\bracket{1, -1, -1},\quad R_{N}^d=\begin{bmatrix}
     -1&     0&     0\\
    0&     0&     -1\\
     0&     -1&     0\end{bmatrix},\\
\Pi_0=\omega_0 J R^{T}_0e_2,\quad\Pi_N^d=\omega_0JR^{d,T}_Ne_2.
\end{gather*}
\end{enumerate}
The optimized performance index and the violation of the constraints
are given in Table \ref{tab:res}. Figures \ref{fig:sati} and
\ref{fig:satii} show the attitude maneuver of the spacecraft
(clockwise direction), the control inputs, the angular
velocity response, and the violation of the terminal boundary
condition according to the number of iterations.

\subsection{Numerical properties}
The neighboring extremal method or the shooting method are
numerically efficient in the sense that the number of optimization
parameters is minimized. But, this approach may tend to have
numerical ill-conditioning~\cite{bk:betts}. A small change in the
initial multiplier can cause highly nonlinear behavior of the
terminal attitude and angular momentum. It is difficult to compute
the Jacobian matrix for Newton iterations accurately, and
consequently, the numerical error may not converge.

However, the numerical examples presented in this section show
excellent numerical convergence properties. They exhibit a quadratic
rate of convergence. This is because the proposed computational
algorithms on $\SO$ are geometrically exact and numerically
accurate. The attitude dynamics of a rigid body arises from
Hamiltonian mechanics, which have neutral stability. The adjoint
system is also neutrally stable. The proposed Lie group variational
integrators and the discrete multiplier equations, obtained from
variations expressed in the Lie algebra, can preserve the neutrally
stability property. Therefore the sensitivity derivatives are
computed accurately.

\section{Conclusions}
A discrete optimal control problem for the attitude dynamics of a
rigid body in the presence of an attitude dependent potential is
studied. The performance index is the $l_2$ norm of external control
inputs and boundary conditions on the attitude and the angular
momentum are prescribed. The attitude is represented by a rotation
matrix in the Lie group, $\SO$. This paper proposes three levels of
geometrically exact computations on $\SO$ to solve the optimal
control problem; Lie group variational integrator, discrete-time
necessary conditions for optimality, and discrete-time sensitivity
derivatives.

The Lie group variational integrator obtained from a discrete
variational principle preserves the geometric features of the
attitude dynamics of the rigid body. It exhibits symplectic and
momentum preservation properties, and good energy behavior
characteristic of variational integrators. Since the rotation
matrices are updated by a group operation, the Lie group structure
is also preserved.

The necessary conditions of optimality are derived by a variational
principle. The Lie group variational integrators are imposed as
constraints, and the variation of the rotation matrices are
expressed in terms of Lie algebra elements. The proposed discrete
optimality conditions are the basis for a numerically efficient
computational algorithms for the optimal attitude control problem,
since the implicit part of the optimality conditions occurs in a
single equation of one variable. This implicit condition can be
solved easily by Netwon iteration. Other algorithms require
iteration on the entire discrete time trajectory simultaneously.

The necessary conditions are expressed as a two point boundary value
problem on $\T^*\SO$ and its dual space. The sensitivity derivatives
are developed in the Lie algebra, and the two point boundary value
problem is solved using a neighboring extremal method. The
neighboring extremal method is efficient for this class of optimal
control problems because the resulting problem of satisfying the
terminal boundary conditions has a small number of variables. The
main disadvantage is that a small change in the initial multipliers
can produce a very large change in the terminal condition. This can
result in numerical ill-conditioning. The nonlinearity also makes it
hard to construct an accurate estimate of the Jacobian matrix that
is needed for a Newton iteration. In this paper, the two point
boundary problem is solved efficiently. The error in the terminal
boundary conditions converges exponentially to machine precision.
This is because the sensitivity derivatives are computed accurately
in the Lie algebra of $\SO$.

Numerical results for an optimal attitude control problem involving
an underactuated axially symmetric 3D pendulum and for an optimal
attitude control problem involving a fully actuated spacecraft on a
circular orbit are given. The boundary conditions are chosen such
that the resulting maneuvers are large angle attitude maneuver. It
is shown that the proposed numerical computations on $\SO$ are
geometrically exact and highly efficient.

\newpage
\bibliography{opt}
\bibliographystyle{jota}

\newpage
\section*{Appendix} The linearized equations of motion for the
multiplier equations \refeqn{updatelam} are given by
\begin{align}
\begin{bmatrix} \delta \lambda_{k}^1 \\ \delta\lambda_{k}^2 \end{bmatrix}
= \begin{bmatrix} \delta\mathcal{A}_{k+1}^T & \delta\mathcal{C}_{k+1}^T \\
\delta\mathcal{B}_{k+1}^T & \delta\mathcal{D}_{k+1}^T\end{bmatrix}
\begin{bmatrix} \lambda_{k+1}^1 \\ \lambda_{k+1}^2\end{bmatrix}
+ \begin{bmatrix} \mathcal{A}_{k+1}^T & \mathcal{C}_{k+1}^T \\
\mathcal{B}_{k+1}^T & \mathcal{D}_{k+1}^T\end{bmatrix}
\begin{bmatrix} \delta\lambda_{k+1}^1 \\
\delta\lambda_{k+1}^2\end{bmatrix}.\label{eqn:delcs}
\end{align}
In this appendix, we derive expressions for these variations, and we
summarize the results. Here we repeatedly use the properties,
$S(x)y=-S(y)x$, $S(Rx)=RS(x)R^T$ for $x,y\in\Re^3$, $R\in\SO$.

\paragraph{Variation $\delta\mathcal{A}_{k}^T\lambda^1_k$:}Using \refeqn{xik}, the variation
$\delta\mathcal{A}_{k}^T\lambda^1_k$ is given by
\begin{align}
\delta\mathcal{A}_{k}^T\lambda^1_k & = \delta F_k\lambda^1_k,\nonumber\\
& = -F_k S(\lambda^1_k) \mathcal{B}_k \delta\Pi_k.\label{eqn:delAAk}
\end{align}

\paragraph{Variation $\delta\mathcal{C}_{k}^T\lambda^1_k$:}

Since $\mathcal{M}_{k+1}$ depends on $R_{k+1}$,
$\delta\mathcal{M}_{k+1}^T x$ for any $x\in\Re^3$ can be written as
\begin{align}
\delta\mathcal{M}_{k+1}^T x & = \mathcal{N}_{k+1}^T(x) \zeta_{k+1},\nonumber\\
& = \mathcal{N}_{k+1}^T(x) \mathcal{A}_k \zeta_k +
\mathcal{N}_{k+1}^T(x) \mathcal{B}_k\delta\Pi_k,\label{eqn:delMMkpx}
\end{align}
where $\mathcal{N}_{k}(x)\in\Re^{3\times 3}$. Using \refeqn{xik} and
\refeqn{delMMkpx}, the variation
$\delta\mathcal{C}_{k}^T\lambda^1_k$ is given by
\begin{align}
\delta\mathcal{C}_{k}^T\lambda^2_k
& = h\delta F_k\mathcal{M}_{k+1}^T\lambda^2_k + hF_k\delta \mathcal{M}_{k+1}^T\lambda^2_k,\nonumber\\
& = hF_k\mathcal{N}_{k+1}^T(\lambda_k^2) \mathcal{A}_k \zeta_k
+hF_k\braces{-S(\mathcal{M}_{k+1}^T\lambda^2_k)\mathcal{B}_k+\mathcal{N}_{k+1}^T(\lambda_k^2)
\mathcal{B}_k}\delta\Pi_k.\label{eqn:delCCk}
\end{align}

\paragraph{Variation $\delta\mathcal{B}_{k}^T\lambda^2_k$:}To obtain the variation $\delta\mathcal{B}_{k}^T\lambda^1_k$, we rewrite the expression of $\mathcal{B}^T_k\lambda^1_k$ as follows. Using the definition of $\mathcal{B}_k$ in
\refeqn{BBk}, we obtain
\begin{gather*}
\braces{\tr{F_kJ_d}I_{3\times 3}-F_kJ_d}^{T}\mathcal{B}_k^T
\lambda^1_k = hF_k \lambda^1_k.
\end{gather*}
Since $S(x)A^T+AS(x)=S(\braces{\tr{A}I_{3\times 3}-A}^Tx)$ for all
$x\in\Re^3,A\in\Re^{3\times 3}$, the above equation can be written
in matrix form;
\begin{gather}
S(\mathcal{B}_k^T \lambda^1_k)J_dF_k^T+F_kJ_dS(\mathcal{B}_k^T
\lambda^1_k) = S(hF_k \lambda^1_k).\label{eqn:delBlam10}
\end{gather}
Taking a variation of the left hand side of \refeqn{delBlam10},
\begin{align}
S(&\delta(\mathcal{B}_k^T\lambda^1_k))J_dF_k^T+F_kJ_dS(\delta(\mathcal{B}_k^T \lambda^1_k))
+S(\mathcal{B}_k^T \lambda^1_k)J_dS(\mathcal{B}_k\delta\Pi_k)^TF_k^T+F_kS(\mathcal{B}_k\delta\Pi_k)J_dS(\mathcal{B}_k^T \lambda^1_k)\nonumber\\
&=S(\braces{\tr{F_kJ_d}I_{3\times3}-F_kJ_d}^{T}\delta(\mathcal{B}_k^T \lambda^1_k))\nonumber\\
&\quad+S(\braces{\tr{F_kJ_dS(\mathcal{B}_k^T \lambda^1_k)}I_{3\times3}-F_kJ_dS(\mathcal{B}_k^T\lambda^1_k)}F_k\mathcal{B}_k\delta\Pi_k).\label{eqn:43l}
\end{align}
Taking a variation of the right hand side of \refeqn{delBlam10},
\begin{align}
\delta S(hF_k \lambda^1_k) & = S(hF_kS(\mathcal{B}_k\delta\Pi_k)\lambda^1_k+hF_k\delta \lambda^1_k),\nonumber\\
&=S(-hF_kS(\lambda^1_k)\mathcal{B}_k\delta\Pi_k+hF_k\delta\lambda^1_k).\label{eqn:43r}
\end{align}
Using \refeqn{43l} and \refeqn{43r}, $\delta(\mathcal{B}_k\lambda^1_k)$ is given by
\begin{align}
\delta(\mathcal{B}_k^T\lambda^1_k) & = \mathcal{E}_k^T(\lambda_k^1)
\delta\Pi_k + \mathcal{B}_k^T\delta\lambda^1_k,\label{eqn:delBBk}
\end{align}
where $\mathcal{E}_k^T(\cdot):\Re^3\mapsto\Re^{3\times 3}$ is
defined for $x\in\Re^3$ as
\begin{align*}
\mathcal{E}_k^T(x) & =-\braces{\tr{F_kJ_d}I_{3\times 3}-F_kJ_d}^{-T}\\
&\quad \times \bracket{\braces{\tr{F_kJ_dS(\mathcal{B}_k^Tx)}I_{3\times 3}
-F_kJ_dS(\mathcal{B}_k^Tx)}F_k\mathcal{B}_k+hF_kS(x)\mathcal{B}_k}.
\end{align*}

\paragraph{Variation $\delta\mathcal{D}_{k}^T\lambda^2_k$:}The variation $\delta\mathcal{D}_k^T\lambda_k^2$ is given by
\begin{align}
\delta\mathcal{D}_k^T\lambda_k^2 & = \delta F_k\lambda_k^2 + \delta(\mathcal{B}_k^T\braces{-S(F_k^T\Pi_k)+h\mathcal{M}_{k+1}^T})\lambda_k^2,\nonumber\\
& = h\mathcal{B}_k^T\mathcal{N}_{k+1}^T(\lambda_k^2)\mathcal{A}_k
\zeta_k+\mathcal{F}_k^T\delta\Pi_k,\label{eqn:delDDk}
\end{align}
where $\mathcal{F}_k^T\in\Re^{3\times 3}$ is defined as
\begin{align*}
\mathcal{F}_k^T&=-F_kS(\lambda_2^k)\mathcal{B}_k+\mathcal{E}_k^T(\braces{-S(F_k^T\Pi_k)+h\mathcal{M}_{k+1}^T}\lambda_k^2)\\
&\quad+\mathcal{B}_k^T\braces{S(\lambda_k^2)(S(F_k^T\Pi_k)\mathcal{B}_k+F_k^T)
+h\mathcal{N}_{k+1}^T(\lambda_k^2)\mathcal{B}_k}.
\end{align*}

\paragraph{Summary:}Equations \refeqn{delAAk}, \refeqn{delBBk}, \refeqn{delCCk} and
\refeqn{delDDk} are expressions of the variations in \refeqn{delcs}.
Then, \refeqn{delx} and \refeqn{delcs} can be written as
\begin{align*}
x_{k+1} & = A^{11}_k x_k + A^{12}_k \delta\lambda_k,\\
\delta\lambda_k & =A^{21}_{k+1} x_{k+1}+\parenth{A^{11}_{k+1}}^T
\delta\lambda_{k+1},
\end{align*}
where $x_k=[\zeta_k^T,\delta\Pi_k^T]^T$,
$\delta\lambda_k=[\delta\lambda_k^{1,T},\delta\lambda_k^{2,T}]^T\in\Re^6$, and
\begin{align*}
A^{11}_k & = \begin{bmatrix} \mathcal{A}_k & \mathcal{B}_k\\\mathcal{C}_k & \mathcal{D}_k\end{bmatrix},\\
A^{12}_k & = \begin{bmatrix} 0 & 0\\0 & -hBB^T\end{bmatrix},\\
A^{21}_k & = \begin{bmatrix}
hF_k\mathcal{N}_{k+1}^T(\lambda_k^2) \mathcal{A}_k & -F_k S(\lambda^1_k) \mathcal{B}_k+hF_k\braces{-S(\mathcal{M}_{k+1}^T\lambda^2_k)\mathcal{B}_k+\mathcal{N}_{k+1}^T(\lambda_k^2) \mathcal{B}_k} \\
h\mathcal{B}_k^T\mathcal{N}_{k+1}^T(\lambda_k^2)\mathcal{A}_k &
\mathcal{E}_k^T(\lambda_k^1)+\mathcal{F}_k^T
\end{bmatrix}.
\end{align*}
In summary, the linear discrete equations of motion can be written as
\begin{align}
\begin{bmatrix}x_{k+1}\\\delta\lambda_{k+1}\end{bmatrix}
= \begin{bmatrix} A^{11}_k & A^{12}_k\\
-\parenth{A^{11}_{k+1}}^{-T}A^{21}_{k+1}A^{11}_{k} &
\parenth{A^{11}_{k+1}}^{-T}\parenth{I_{3\times
3}-A^{21}_{k+1}A^{12}_{k}}\end{bmatrix}
\begin{bmatrix}x_{k}\\\delta\lambda_{k}\end{bmatrix}.
\end{align}

\newpage\listoftables

\clearpage\newpage
\vspace*{3.5cm}
\renewcommand{\theenumi}{\arabic{enumi}}
\renewcommand{\labelenumi}{\theenumi:}
\newcommand{\tab}{\hspace*{0.8cm}}
\begin{table}[h]
\caption[Newton-Armijo iteration procedures]{}\label{tab:iter}
\hrule\vspace*{0.3cm}
\begin{enumerate}
\item Guess an initial multiplier $\lambda_0$.
\item Find $\Pi_k,R_k,\lambda_k^{1},\lambda_k^{2}$ for $k=1,2,\cdots,N$ using the initial
conditions
 and \refeqn{updatePi1}--\refeqn{updatelam}.
\item Compute the error in satisfaction of the terminal boundary condition;\\
$\zeta_N^{}=S^{-1}\!\parenth{\mathrm{logm}\!
\parenth{R_N^{T} R_{N}^{d}}}$,
$\delta\Pi_N^{}=\Pi_N^{d}-\Pi_N^{}$.\\
$\mathrm{Error}=\norm{[\zeta_N ;\delta\Pi_N]}$.
\item Set $\mathrm{Error}^t=\mathrm{Error},\;\; i=1$.
\item \textbf{while} $\mathrm{Error} > \epsilon_S$.
\item \tab Find a line search direction; $D=\Phi_{12}^{-1}$.
\item \tab Set $c=1$.
\item \tab\tab \textbf{while} $\mathrm{Error}^t > (1-2\alpha c)\mathrm{Error}$
\item \tab\tab Choose a trial initial condition $\lambda_0^{t}=\lambda_0^{}+c D
[\zeta_N;\delta\Pi_N]$.
\item \tab\tab Find $\Pi_k^{},R_k^{},\lambda_k^{1},\lambda_k^{2}$ for $k=1,2,\cdots,N$ using the trial initial
conditions and \refeqn{updatePi1}--\refeqn{updatelam}.
\item \tab\tab Compute the error in satisfaction of the terminal boundary condition\\
\tab\tab$\zeta_N^{t}=S^{-1}\!\parenth{\mathrm{logm}\!\parenth{R_N^{T}
R_{N}^{d}}}$,
$\delta\Pi_N^{t}=\Pi_N^{d}-\Pi_N^{}$.\\
\tab\tab$\mathrm{Error}^t=\norm{[\zeta_N^t ;\delta\Pi_N^t]}$.
\item \tab\tab Set $c=c/2,\;\; i=i+1$.
\item \tab\tab\textbf{end while}
\item \tab Set $\lambda_0=\lambda_0^t$, $\mathrm{Error}=\mathrm{Error}^t$. (accept the trial)
\item \textbf{end while}
\end{enumerate}
\hrule
\end{table}

\clearpage\newpage
\vspace*{7cm}
\begin{table}[h]
\begin{center}
\caption[Optimized performance index and violation of
constraints]{}\label{tab:res}
\begin{tabular}{cc|ccc}\hline\hline
Case & Model & $\mathcal{J}$ & $\norm{\mathrm{logm}(R_N^{d,T}R_N)}$
& $\norm{\Pi_N^d-\Pi_N}$\\\hline
(i)  & \multirow{2}{22mm}{3D Pendulum} & $1.52$ & $1.77\times 10^{-14}$ & $7.08\times 10^{-15}$\\
(ii) &       & $40.22$ & $2.22\times 10^{-16}$ &
$2.55\times10^{-14}$\\\hline
(iii)& \multirow{2}{22mm}{Spacecraft}  & $23.35$ & $2.90\times 10^{-15}$ & $5.13\times 10^{-15}$\\
(iv) &       & $70.74$ & $7.31\times 10^{-15}$ & $1.48\times 10^{-14}$\\
\hline\hline
\end{tabular}
\end{center}
\end{table}

\clearpage\newpage \setcounter{lofdepth}{2} \listoffigures

\clearpage\newpage\vspace*{1cm}
\renewcommand{\xyWARMinclude}[1]{\includegraphics[width=0.6\textwidth]{#1}}
\begin{figure*}[h]
    \centerline{\subfigure[3D Pendulum][]{
    $$\begin{xy}
    \xyWARMprocessEPS{3dpend}{eps}
    \xyMarkedImport{}
    \xyMarkedMathPoints{1-4}
    \end{xy}$$\label{fig:3dpend}}}
    \centerline{\subfigure[Spacecraft on a circular orbit][]{
    $$\begin{xy}
    \xyWARMprocessEPS{coordinate}{eps}
    \xyMarkedImport{}
    \xyMarkedMathPoints{1-9}
    \end{xy}$$\label{fig:sat}}
    }\caption[Rigid body models]{}\label{fig:body}
\end{figure*}

\clearpage\newpage\vspace*{1cm}

\clearpage\newpage\vspace*{1cm}

\clearpage\newpage\vspace*{1cm}

\clearpage\newpage\vspace*{1cm}

\end{document}